\documentclass[11pt]{amsart}
\usepackage{latexsym, amsmath, amssymb,amsthm,amsopn,amsfonts,mathtools}
\usepackage{version}
\usepackage{epsfig,graphics,color,graphicx,graphpap}
\usepackage{amssymb}
\usepackage{todonotes}
\usepackage{listings}
\usepackage{mathrsfs}
\usepackage[citecolor=blue]{hyperref}
\usepackage{soul,xcolor}

\usepackage[framemethod=tikz]{mdframed}
\newcommand{\redsout}{\bgroup\markoverwith{\textcolor{red}{\rule[0.5ex]{2pt}{.4pt}}}\ULon}

\numberwithin{equation}{section}
\newtheorem{theorem}{Theorem}[section]
\newtheorem*{theorem*}{Theorem}
\newtheorem{corollary}[theorem]{Corollary}
\newtheorem{proposition}[theorem]{Proposition}
\newtheorem*{theorem_0*}{Proposition}
\newtheorem{lemma}[theorem]{Lemma}

\newtheorem{remark}{Remark}[section]

\title[]{Recovery of Coefficients in a Schr\"{o}dinger System with Nonlinear Coupling}
\author[P. Arrepu]{Pranav Arrepu}
\address{Department of Mathematics, University of California Santa Barbara, Santa Barbara, CA 93106, USA}
\email{pranav\_arrepu@ucsb.edu}

\author[H. Zhou]{Hanming Zhou}
\address{Department of Mathematics, University of California Santa Barbara, Santa Barbara, CA 93106, USA}
\email{hzhou@math.ucsb.edu}

\begin{document}

\begin{abstract}
We consider the inverse problem of recovering time-dependent coefficients in a system of Schr\"{o}dinger equations with coupled nonlinear terms. We first address the well-posedness of the system for small boundary data and trivial initial data. Then, we show that the linearized Neumann data stably determines the nonlinear coefficients provided that, at any given time, they are known in some neighborhood of the boundary. The proof relies on linearization and complex geometric optics (CGO) solutions for the linear Schr\"{o}dinger equation. We conclude by outlining an approach based on the direct use of Carleman estimates to recover time-independent coefficients in the setting of trivial boundary data but nontrivial initial data.
\end{abstract}

\maketitle

\section{Introduction} \label{introduction}
Letting $\Omega$ be a smooth bounded domain in $\mathbb{R}^{n}\;(n = 2, 3)$ and denoting the linear Schr\"{o}dinger operator by $\mathcal{T}_{j} \coloneqq i\partial_t + \Delta + p_{j}(x,t)$, we seek to study the following system of coupled nonlinear Schr\"{o}dinger equations:
\begin{equation} \label{two_component_GPEqn}
\begin{cases}
\mathcal{T}_{1} u_{1} + q_{1}(x,t)\left|u_{1}\right|^{2}u_{1} + \tilde{q}_{1}(x, t)|u_{2}|^{2}u_{1} = 0 \text{ on } \Omega \times (0, T),\\
\mathcal{T}_{2} u_{2} + q_{2}(x,t)\left|u_{2}\right|^{2}u_{2} + \tilde{q}_{2}(x, t)|u_{1}|^{2}u_{2} = 0 \text{ on } \Omega \times (0, T),\\
u_{1}(x,t) = g_{1}(x,t),\,u_{2}(x,t) = g_{2}(x,t) \text{ on } \partial \Omega \times (0,T),\\
u_{1}(x,0) = u_{2}(x,0) = 0 \text{ on } \Omega \times \{0\}.
\end{cases}
\end{equation}
When the coefficients $q_{j}, \tilde{q}_{j}\;(j = 1,2)$ are real constants, it is noted in \cite{SK2020} that \eqref{two_component_GPEqn} could model a two-component Bose-Einstein condensate with wavefunction $\mathbf{u} = (u_{1},u_{2})$ or coupled modes in birefringent media, where $u_{1}$ and $u_{2}$ are the electric field intensities along each birefringence axis. It is also noted in \cite{LK2013} that variable coefficients can arise in the nonlinear Schr\"{o}dinger equation when modelling optical pulses through nonlinear optical fibers.\\ \indent
We begin by recording some notation used in the paper. We will denote $Q \coloneq \Omega \times (0, T)$ and $\Sigma \coloneq \partial \Omega \times (0, T)$. Taking $\Gamma_{j}\,(j = 1,2)$ arbitrary relatively open subsets of $\partial \Omega$, we also denote $\Sigma_{j} \coloneq \Gamma_{j} \times (0, T)$. Throughout the paper, we will use the Sobolev space
\begin{equation} \label{time_space_convention}
H^{r,s}(S \times (0, T)) \coloneqq H^{r}(0, T; L^{2}(S)) \cap L^{2}(0, T;H^{s}(S)).
\end{equation}
For a given Banach space $B$, we will consider when needed the product space $B^{l}$ for $l \in \mathbb{N}$ and equip it with the norm given by
\begin{equation} \label{product_space_norm_0}
\left \Vert \left(b_{1},\ldots,b_{l}\right) \right \Vert_{B^{l}} \coloneqq \sum_{i = 1}^{l}\left \Vert b_{i} \right \Vert_{B}.
\end{equation}
We will also use the permutation $\iota:\{1,2\} \longrightarrow \{1,2\}$ defined by $\iota(1) = 2$ and $\iota(2) = 1$. \\ \indent
We are interested in the following inverse problem:\\
\textbf{(IP)} Do the coefficients $q_{j}, \tilde{q}_{j}\; (j = 1,2)$ depend continuously on the Neumann data
\begin{equation} \label{main_inverse_stability_measurement}
\left(\partial u_{1}/\partial \nu\big|_{\Sigma_{1}},\partial u_{2}/\partial \nu\big|_{\Sigma_{2}}\right)?
\end{equation}
That is, we seek to stably recover the nonlinear coefficients of \eqref{two_component_GPEqn} from the Neumann data \eqref{main_inverse_stability_measurement}. This question is related to the Calder\'{o}n problem for the conductivity equation $\nabla \cdot (\gamma\nabla u) = 0$, which seeks to determine the coefficient $\gamma$ from the knowledge of the Dirichlet-to-Neumann map (DN-map) $f \longmapsto \gamma\left(\partial u/\partial \nu \right)\big|_{\partial\Omega}$, where $f$ is the Dirichlet data. In \cite{SU1987}, the unique determination of $\gamma$ from the DN-map was addressed in the affirmative for a smooth $\gamma$ with positive lower bound and since then, a large body of work has emerged to address the Calder\'{o}n problem for other equations. We briefly mention references for the dynamical (time-dependent) Schr\"{o}dinger equation and coupled PDE systems.\\ \indent
The proof of the uniqueness result in \cite{SU1987} relies on a transformation of the conductivity equation into a time-independent Schr\"{o}dinger equation. However, the Calder\'{o}n problem for the dynamical linear Schr\"{o}dinger equation began with the work of \cite{BP2002}. In this work, the time-independent electromagnetic potential was stably recovered by the Neumann data measured on a part of the boundary satisfying a geometric constraint. Subsequent work on the recovery of coefficients in the linear dynamical Schr\"{o}dinger equation includes \cite{Eskin2008, MOR2008, BC2009, YY2010, KSS2014, CKS2015, BKS2016, HKSY2019, KS2019, BB2020}. Inverse problems for the dynamical Schr\"{o}dinger equation with polynomial nonlinearities began with the work of \cite{LOST2022}, where the time-dependent linear and nonlinear coefficients were recovered from the source-to-solution map. Polynomial nonlinearities were also considered in \cite{LLZ2023,LUY2024, KLV2026}. In \cite{LLZ2023}, the authors considered the recovery of quadratic and Gross-Pitaevskii nonlinear coefficients from the DN-map. The authors of \cite{LUY2024, KLV2026} treat a nonlinear magnetic Schr\"{o}dinger equation and recover both the linear (including the magnetic and electric potentials) and nonlinear coefficients from the DN-map. The nonlinearity considered in \cite{LUY2024} is a finite sum of polynomials in the solution variable and its conjugate while the nonlinearity considered in \cite{KLV2026} is an infinite sum of polynomials in the solution variable. The results in \cite{LOST2022, LLZ2023, LUY2024, KLV2026} apply to the setting of time-dependent coefficients, nontrivial boundary data, and trivial initial data. In \cite{AZ2025}, a modified setting (with trivial boundary data but nontrivial initial data) was considered for the Schr\"{o}dinger equation with a nonlinearity that was a locally analytic function of the solution variable and its conjugate. In \cite{AZ2025}, the time-independent nonlinear coefficient was recovered from the linearized Neumann data. We mention that while the studies \cite{LLZ2023, LUY2024, AZ2025, KLV2026} only require measurements of the Neumann data on an arbitrary relatively open subset of the boundary, they require the coefficients to be known \textit{a priori} in a neighborhood of the entire boundary. \\ \indent
Previous work on the Calder\'{o}n problem for PDE systems with coupling through linear terms includes \cite{LT2011, WY2017, KRS2018, DY2019, DZ2023, HKS2025, Liu2025, SZZ2025, BL2026, BV2026, Saci2026_1, Saci2026_2}. We also mention \cite{LJ2026} which concerns an inverse problem for a single nonlinear elliptic equation but whose analysis relies on boundary determination and unique continuation for a coupled linear elliptic system. The studies in \cite{DY2019, DZ2023, KRS2018, HKS2025, Saci2026_1, Saci2026_2} are concerned with systems of linear Schr\"{o}dinger equations and they establish the stable recovery of coefficients through measurements of the solution in a subset of the interior of the domain or by measurements of the Neumann data on a subset of the boundary. In either case, the proof relies on an appropriate Carleman estimate requiring geometric conditions on the set over which measurements take place. The use of Carleman estimates in the style of \cite{DY2019, DZ2023, KRS2018, HKS2025, Saci2026_1, Saci2026_2} (and for example \cite{BP2002, MOR2008, BC2009, YY2010, HKSY2019, AZ2025, SR2023} in the case of a single equation) goes back to A. L. Bukhgeim and M. V. Klibanov in \cite{BK1981}. \\ \indent
Studies on the Calder\'{o}n problem for PDE systems with coupling through nonlinear terms include \cite{Waters2019, DLZ2023, DLL2024, KLY2024, CLY2026}. The works \cite{ME2024, QGZLLD2024} use neural networks to numerically reconstruct the coefficients in reaction-diffusion and Schr\"{o}dinger systems (respectively) with nonlinear coupling. To our knowledge, the present paper is the first to study the stable recovery of coefficients in a system of Schr\"{o}dinger equations with coupled nonlinear terms.\\ \indent
We conclude our literature review by noting that the strategy used in \cite{LOST2022, LLZ2023, DLL2024, LUY2024, AZ2025} to study the inverse problem for nonlinear equations is to consider derivatives of the solution with respect to small parameters in the source, boundary, or initial data. The differentiated solution satisfies a linear equation for which the inverse problem is more tractable. This linearization technique goes back to V. Isakov in \cite{Isakov1993} and was further developed in \cite{KLU2018, KU2020, FO2020, LLLS2021}. \\ \indent
In order to address \textbf{(IP)}, we introduce small parameters $\epsilon_{1},\epsilon_{2},\epsilon_{3} \in \mathbb{R}$ into the boundary data of \eqref{two_component_GPEqn} to consider
\begin{equation} \label{two_component_GPEqn_small_data}
\begin{cases}
\mathcal{T}_{j} u_{j,k} + q_{j,k}(x,t)\left|u_{j,k}\right|^{2}u_{j,k} + \tilde{q}_{j,k}(x, t)|u_{\iota(j),k}|^{2}u_{j,k} = 0 \text{ on } Q,\\
u_{j,k}(x,t) = \left(\epsilon_{1}f^{(j)}_{1} + \epsilon_{2}f^{(j)}_{2} + \epsilon_{3}f^{(j)}_{3}\right)(x,t) \text{ on } \Sigma,\\
u_{j,k}(x,0) = 0 \text{ on } \Omega \times \{0\},
\end{cases}
\end{equation}
where $j \in \{1,2\},\, k \in \{1,2\},$ and
\begin{equation} \label{boundary_neighborhood_0}
f^{(j)}_{l} \in \mathcal{H}^{9/2,9/2}(\Sigma) \coloneqq \left\{g \in H^{9/2,9/2}(\Sigma) \,\big|\, \partial^{m}_{t}g(\cdot,0) = 0 \text{ for integers } 0 \leq m < 9/2\right\}
\end{equation}
for $l \in \{1,2,3\}$ are real-valued. For \eqref{two_component_GPEqn_small_data}, we assume that $p_{j} \in W^{11,\infty}(Q), q_{j,k} \in H^{3}(Q), \tilde{q}_{j,k} \in H^{3}(Q)$ are real-valued and that there exist constants $N_{j}, M_{j} > 0$ such that
\begin{equation} \label{coefficient_smallness_assumption}
\begin{split}
\Vert p_{j} \Vert_{W^{11,\infty}(Q)} \leq N_{j},\\
\Vert q_{j,k} \Vert_{H^{3}(Q)}, \Vert \tilde{q}_{j,k} \Vert_{H^{3}(Q)} \leq M_{j}.
\end{split}
\end{equation}
We also suppose that there exists a relatively open neighborhood $\mathcal{O}^{(j)} \subset \Omega$ of $\partial \Omega$ such that on $\mathcal{O}^{(j)} \times (0, T)$, we have 
\begin{equation} \label{near_the_boundary_a_priori_assumption}
\begin{split}
p_{j}(x,t) = p_{j}(x),\\
q_{j,1} \equiv q_{j,2},\\
\tilde{q}_{j,1} \equiv \tilde{q}_{j,2}.
\end{split}
\end{equation}
That is, the coefficient $p_{j}$ is time-independent on $\mathcal{O}^{(j)} \times (0, T)$. Then, our main result is as follows.
\begin{theorem} \label{main_inverse_stability_theorem}
Let $\Omega$ be a smooth bounded domain in $\mathbb{R}^{n}\,(n = 2, 3)$ and let $\mathcal{O}^{(j)} \subset \Omega$ be a relatively open neighborhood of $\partial\Omega$ such that $\partial\mathcal{O}^{(j)} \setminus \partial\Omega$ is $C^{2}$. For \eqref{two_component_GPEqn_small_data}, suppose $p_{j} \in W^{11,\infty}(Q), q_{j,k} \in H^{3}(Q), \tilde{q}_{j,k} \in H^{3}(Q)$ are real-valued and that the coefficients satisfy \eqref{coefficient_smallness_assumption}-\eqref{near_the_boundary_a_priori_assumption}. Suppose that the boundary data of \eqref{two_component_GPEqn_small_data} are real-valued and satisfy \eqref{boundary_neighborhood_0}.\\ \indent Then, for $j \in \{1,2\}$, $l \in \{1,2,3\}$, and any $T^{*} \in (0, T)$, there exists $\delta_{*,j} \in (0, 1)$ such that if \eqref{two_component_GPEqn_small_data} is well-posed for all sufficiently small $\epsilon_{l} \in \mathbb{R}$, and
\begin{equation} \label{Neumann_data_measurement_0}
\delta_{j,\Gamma_{j}} \coloneqq \left\Vert \partial^{3}_{\epsilon_{1}\epsilon_{2}\epsilon_{3}}\left(\partial_{\nu}\left(u_{j,1} - u_{j,2}\right)\right)\Big|_{\epsilon_{1} = \epsilon_{2} = \epsilon_{3} = 0}\right\Vert_{L^{2}(\Sigma_{j})} \in (0,\delta_{*,j}),
\end{equation}
we have for some constants $\alpha_{1},\alpha_{2}, \alpha_{3} > 0$ that
\begin{equation} \label{main_stability_estimate_statement}
\begin{split}
\left \Vert q_{j,1} - q_{j,2} \right \Vert^{2}_{L^{2}(\Omega\setminus\mathcal{O}^{(j)} \times (0, T^{*}))}
\leq C_{j}\left[\left(\ln \delta_{j,\Gamma_{j}}^{-\alpha_{1}} \right)^{-\alpha_{2}} + \delta_{j,\Gamma_{j}}^{\alpha_{3}}\right],\\
\left \Vert \tilde{q}_{j,1} - \tilde{q}_{j,2} \right \Vert^{2}_{L^{2}(\Omega\setminus\mathcal{O}^{(j)} \times (0, T^{*}))}
\leq \tilde{C}_{j}\left[\left(\ln \delta_{j,\Gamma_{j}}^{-\alpha_{1}} \right)^{-\alpha_{2}} + \delta_{j,\Gamma_{j}}^{\alpha_{3}}\right],
\end{split}
\end{equation}
where $C_{j}, \tilde{C}_{j} > 0$ are constants that depend on $\Omega, \Gamma_{j}, \mathcal{O}^{(j)}, T, T^{*}, N_{j}, M_{j}$.
\end{theorem}
The proof of Theorem $1.1$ relies on the linearization of \eqref{two_component_GPEqn_small_data} with respect to the small parameters $\epsilon_{1}, \epsilon_{2}, \epsilon_{3} \in \mathbb{R}$. In particular, we consider the first-order linearization
\begin{equation} \label{heuristic_first_order_linearization}
v^{(l)}_{j,k} \coloneqq \partial_{\epsilon_{l}}u_{j,k}\big|_{\epsilon_{l} = 0}
\end{equation}
and the third-order linearization
\begin{equation} \label{heuristic_kth_order_linearization}
w^{(1,2,3)}_{j,k} \coloneqq \partial^{3}_{\epsilon_{1}\epsilon_{2}\epsilon_{3}}u_{j,k}\big|_{\epsilon_{1}=\epsilon_{2}=\epsilon_{3}=0}.
\end{equation}
The partial derivatives in \eqref{heuristic_first_order_linearization}-\eqref{heuristic_kth_order_linearization} are taken in the sense of the mapping from the boundary data of \eqref{two_component_GPEqn_small_data} to the unique solution of \eqref{two_component_GPEqn_small_data} (the existence of these derivatives will be justified through Corollary \ref{nonlinear_solution_Taylor_series_corollary}).  From here, we apply the method developed in \cite{LLZ2023} to recover nonlinear coefficients. In particular, we begin by multiplying the function
\begin{equation}
w^{(1,2,3)}_{j} \coloneqq w^{(1,2,3)}_{j,1} - w^{(1,2,3)}_{j,2}
\end{equation}
by a smooth cut-off function $\chi(x)$ that restricts $w^{(1,2,3)}_{j}$ to $\mathcal{O}^{(j)} \times (0, T)$. The purpose of this step is to later apply a unique continuation estimate from \cite{BC2009} that relates the $H^{1}$-norm of $w^{(1,2,3)}_{j}(\cdot, t)$ on a subset of $\mathcal{O}^{(j)}$ to the linearized Neumann data defined by \eqref{Neumann_data_measurement_0}. Then, with $v_{0,j}$ denoting a solution to the adjoint problem for the first-order linearization, we use the equation satisfied by $\chi w^{(1,2,3)}_{j}$ to derive an integral identity involving the products of $q_{1} - q_{2}, \tilde{q}_{1} - \tilde{q}_{2}$ with the functions $\left\{v^{(l)}_{j,k},\overline{v_{0,j}}\right\}$ and the product of $[\Delta, \chi]w^{(1,2,3)}_{j}$ with $\overline{v_{0,j}}$. We next substitute complex geometric optics (CGO) solutions constructed for the linear Schr\"{o}dinger equation in \cite{LLZ2023} for the functions $v^{(l)}_{j,k}$ and $v_{0,j}$ appearing in the integral identity. By carefully picking the vectors defining the exponential phases in the CGO solutions, we can derive an estimate of the Fourier transform of either $q_{1} - q_{2}$ or $\tilde{q}_{1} - \tilde{q}_{2}$ in terms of an integral involving $[\Delta, \chi]w^{(1,2,3)}_{j}$. This step uses an observation from \cite{LUY2024} that allows one to separately recover the coefficients of polynomials of the same order. Then, we use the previously mentioned unique continuation estimate to introduce the linearized Neumann data defined by \eqref{Neumann_data_measurement_0} through the $[\Delta, \chi]w^{(1,2,3)}_{j}$ term.
\begin{remark} \label{linear_nonlinear_coefficient_assumptions_0}
\hfill
\begin{itemize}
\item[(a)] Our requirement that $q_{j,k},\tilde{q}_{j,k} \in H^{3}(Q)$ is related to the fact that $m = 3$ is the smallest integer for which $H^{m}(Q)$ is a Banach algebra when $\Omega$ is a smooth bounded domain in $\mathbb{R}^{2}$ or $\mathbb{R}^{3}$\footnote{See Theorem $4.39$ in \cite{AF2003}.}. For example, we will use the Banach algebra property in the proof of the well-posedness of \eqref{two_component_GPEqn} for small boundary data. The assumption $p_{j} \in W^{11,\infty}(Q)$ will later ensure that the CGO solutions we use for the linear Schr\"{o}dinger equation have sufficient regularity.
\item[(b)] The restriction that $p_{j}$ is real-valued is used in the proof of an integral identity stated in Proposition \ref{main_integral_identity}. Moreover, we stipulate $p_{j}(x,t) = p_{j}(x)$ on $\mathcal{O}^{(j)} \times (0, T)$ in order to use a unique continuation result stated in Proposition \ref{unique_continuation_result_BC2009}. The assumption that $p_{j}, q_{j,k}, \tilde{q}_{j,k}, f^{(j)}_{l}$ are real-valued along with the fact that $\partial\mathcal{O}^{(j)}\setminus\partial\Omega$ is $C^{2}$ are used in the proof of this unique continuation result. Further details about the unique continuation result and its proof are provided in Remark \ref{BC2009_unique_continuation_modifications}-\ref{real_valued_data_assumptions}.
\item[(c)] One can uniquely determine the coefficient $p_{j}$ by applying existing results in the literature to the first-order linearization defined by \eqref{heuristic_first_order_linearization}. This is discussed in Remark \ref{linear_coefficient_recovery}. A unique determination result for the coefficients $q_{j,k}, \tilde{q}_{j,k}$ is discussed in Remark \ref{unique_determination_result_0}.
\item[(d)] The effect of the nonlinear coupling in \eqref{two_component_GPEqn_small_data} is present in the quantity $\delta_{j,\Gamma_{j}}$ defined by \eqref{Neumann_data_measurement_0}. Indeed, we will see in Section \ref{rigorous_linearization} that the function defined by \eqref{heuristic_kth_order_linearization} satisfies an equation (namely \eqref{kth_order_linearization_time-dependent}) whose source term involves derivatives of the nonlinearity
\begin{equation*}
q_{j,k}(x,t)\left|u_{j,k}\right|^{2}u_{j,k} + \tilde{q}_{j,k}(x, t)|u_{\iota(j),k}|^{2}u_{j,k}
\end{equation*}
that appears in \eqref{two_component_GPEqn_small_data}. As a result, $\delta_{j,\Gamma_{j}}$ implicitly depends on both $u_{1,k}$ and $u_{2,k}$.
\end{itemize}
\end{remark}
\begin{remark}
In \eqref{two_component_GPEqn}, we consider nonlinear terms that are third degree polynomials in $\{u_{1}, \overline{u_{1}},u_{2},\overline{u_{2}}\}$. However, the methods in the present paper can be carried over to the case of more general nonlinear terms that are locally analytic in $\{u_{1}, \overline{u_{1}},u_{2},\overline{u_{2}}\}$. We refer the reader to \cite{AZ2025} for a treatment of such nonlinearities for a single equation.
\end{remark}
The rest of the paper is organized as follows. We begin by establishing the well-posedness of \eqref{two_component_GPEqn} for small boundary data in Section \ref{well-posedness_linear_nonlinear}. In Section \ref{rigorous_linearization}, we linearize \eqref{two_component_GPEqn_small_data} with respect to the small parameters $\epsilon_{1},\epsilon_{2},\epsilon_{3} \in \mathbb{R}$. In Section \ref{stability_estimate_proof}, we prove the main result by using the previously discussed method from \cite{LLZ2023} to recover nonlinear coefficients via an integral identity involving the third-order linearization restricted to $\mathcal{O}^{(j)} \times (0, T)$. In Section \ref{Carleman_estimate_recovery_section}, we outline a method from \cite{AZ2025} based on the direct application of Carleman estimates to obtain a stability estimate similar to that of Theorem \ref{main_inverse_stability_theorem} when \eqref{two_component_GPEqn} is equipped with time-independent coefficients, trivial boundary data, and nontrivial initial data.
\bigskip

{\bf Acknowledgments:} The authors are partly supported by the NSF grant DMS-2408369 and Simons Foundation Travel Support for Mathematicians MPS-TSM-00008046.

\section{The Forward Problem} \label{well-posedness_linear_nonlinear}
We seek to show the well-posedness of \eqref{two_component_GPEqn} in $\left(H^{3}(Q)\right)^{2}$ in Theorem \ref{nonlinear_well-posedness}. Our argument will rely on Theorem \ref{linear_well-posedness}, which addresses the well-posedness of the linear Schr\"{o}dinger equation in $H^{m}(Q)$ for integers $m \geq 1$. We begin by noting the following result from Lemma $2.2$ and Remark $2.3$ in \cite{KS2019}.
\begin{lemma} \label{KS2019_Lemma_2.2_statement}
Suppose $p \in W^{1,\infty}(Q)$ and $F \in H^{1}(0, T; L^{2}(\Omega))$ satisfies $F(\cdot, 0)=  0$. Then,
\begin{equation} \label{linear_equation_0_0}
\begin{cases}
(i\partial_{t} + \Delta + p(x,t))v = F(x,t) \text{ on } Q,\\
v(x,t) = 0 \text{ on } \Sigma,\\
v(x,0) = 0 \text{ on } \Omega \times \{0\}.
\end{cases}
\end{equation}
has a unique solution $v \in C\left([0, T]; H^{2}(\Omega) \cap H^{1}_{0}(\Omega)\right)$ satisfying the estimates
\begin{equation} \label{KS2019_Lemma_2.2_estimate}
\begin{split}
\left \Vert v \right \Vert_{L^{\infty}(0, T; L^{2}(\Omega))} \leq C \left \Vert F \right \Vert_{L^{2}(Q)},\\
\left \Vert v \right \Vert_{C([0, T]; H^{2}(\Omega))} + \left \Vert v \right \Vert_{C^{1}([0, T]; L^{2}(\Omega))} \leq C \left \Vert F \right \Vert_{H^{1}(0, T; L^{2}(\Omega))}
\end{split}
\end{equation}
for some constant $C = C\left(\Omega, T, \left\Vert p \right\Vert_{W^{1,\infty}(Q)}\right) > 0$.
\end{lemma}
Now we note a linear well-posedness result that is similar to Proposition $2.1$ from \cite{LLZ2023} and Lemma $4$ from \cite{LOST2022}. We still provide the details of the proof to emphasize that it suffices to take $p \in W^{m,\infty}(Q)$ (whereas the authors of \cite{LLZ2023, LOST2022} assume smoothness for the linear coefficient $p$).
\begin{theorem} \label{linear_well-posedness}
Let $\Omega$ be a smooth bounded domain in $\mathbb{R}^{n}\; (n \geq 2)$. For $m \geq 1$ an integer, suppose that $p \in W^{m,\infty}(Q)$, $h \in H^{m}(Q)$ satisfies $\partial^{k}_{t}h(\cdot, 0) = 0$ for integers $0 \leq k \leq m$, and $g \in H^{m + (3/2),m + (3/2)}(\Sigma)$ satisfies $\partial^{m}_{t}g(\cdot,0) = 0$ for integers $0 \leq k < m + (3/2)$. Then,
\begin{equation} \label{linear_equation_0}
\begin{cases}
(i\partial_{t} + \Delta + p(x,t))v = h(x,t) \text{ on } Q,\\
v(x,t) = g(x,t) \text{ on } \Sigma,\\
v(x,0) = 0 \text{ on } \Omega \times \{0\}.
\end{cases}
\end{equation}
has a unique solution $v \in H^{m}(Q)$ satisfying the estimate
\begin{equation} \label{linear_equation_boundary_data_estimate}
\Vert v \Vert_{H^{m}(Q)} \leq C \left( \Vert g \Vert_{H^{m + (3/2), m + (3/2)}(\Sigma)} + \Vert h \Vert_{H^{m}(Q)}\right)
\end{equation}
for some constant $C = C\left(\Omega, T, \left\Vert p \right\Vert_{W^{m,\infty}(Q)} \right) > 0$. Moreover, we have that $\partial^{k}_{t}v(\cdot, 0) = 0$ for integers $0 \leq k \leq m$.
\end{theorem}
\begin{proof}
Using Theorem $2.3$ in Chapter $4$ of \cite{LM1972}, we have that there exists $v_{2} \in H^{m+2,m+2}(Q)$ such that
\begin{equation} \label{linear_equation_1_0}
\partial^{k}_{t}v_{2}(\cdot,0) = 0
\end{equation}
for integers $0 \leq k < m + (3/2)$ and that $g = v_{2}|_{\Sigma}$. The same theorem in tandem with the Bounded Inverse Theorem (e.g. Theorem $8.34$ in \cite{RR2004}) implies the inequality
\begin{equation} 
\Vert v_{2} \Vert_{H^{m+2,m + 2}(Q)} \leq C\Vert g \Vert_{H^{m + (3/2), m + (3/2)}(\Sigma)}
\end{equation}
for some constant $C = C(\Omega, T) > 0$. We also have from Proposition $2.3$ in Chapter $4$ of \cite{LM1972} the fact that
\begin{equation} \label{v_2_estimate_step_0}
\Vert v_{2} \Vert_{H^{m + 2}(Q)} \leq C\Vert v_{2} \Vert_{H^{m+2,m + 2}(Q)} \left( \leq C\Vert g \Vert_{H^{m + (3/2), m + (3/2)}(\Sigma)}\right)
\end{equation}
for some constant $C = C(\Omega, T) > 0$. Since $p \in W^{m,\infty}(Q)$, we have that $-(i\partial_{t} + \Delta + p)v_{2} \in H^{m}(Q)$. We will now consider the solution $v_{1}(x,t)$ to
\begin{equation} \label{linear_equation_1}
\begin{cases}
(i\partial_{t} + \Delta + p(x,t))v_{1} = F(x,t) \coloneqq h(x,t) - (i\partial_{t} + \Delta + p)v_{2} \text{ on } Q,\\
v_{1}(x,t) = 0 \text{ on } \Sigma,\\
v_{1}(x,0) = 0 \text{ on } \Omega \times \{0\}
\end{cases}
\end{equation}
and show that $v \coloneqq v_{1} + v_{2}$ satisfies the estimate in \eqref{linear_equation_boundary_data_estimate}. \\ \indent
From Lemma \ref{KS2019_Lemma_2.2_statement}, we know that \eqref{linear_equation_1} has a unique solution in $L^{\infty}(0, T; L^{2}(\Omega))$ satisfying the estimate
\begin{equation} \label{linear_equation_1_estimates}
\Vert v_{1}\Vert_{L^{\infty}(0, T; L^{2}(\Omega))} \leq C\Vert F \Vert_{L^{2}(Q)}
\end{equation}
for some constant $C = C\left(\Omega, T, \left\Vert p \right\Vert_{W^{1,\infty}(Q)}\right) > 0$ provided that $F(\cdot,0) = 0$. In fact, recall that we have assumed $\partial^{k}_{t}h(\cdot, 0) = 0$ for integers $0 \leq k \leq m$. So, it follows from \eqref{linear_equation_1_0} that $\partial^{k}_{t}F(\cdot,0) = 0$ for integers $0 \leq k \leq m$. \\ \indent
We will now show that
\begin{equation} \label{v_1_estimate_step_0}
\begin{split}
\left\Vert \partial^{\alpha}_{t}\partial^{\beta}_{x}v_{1} \right\Vert_{L^{2}(Q)} \leq C\left\Vert F \right\Vert_{H^{m}(Q)}\text{ for } 0 \leq \alpha + |\beta| \leq m,\\
\partial^{\alpha}_{t}\partial^{\beta}_{x}v_{1}(\cdot, 0) = 0 \text{ for } 0 \leq \alpha + |\beta| \leq m,
\end{split}
\end{equation}
where $\beta$ is a multi-index and $C = C(\Omega, T, \left\Vert p \right\Vert_{W^{m,\infty}(Q)}) > 0$. Indeed, \eqref{v_1_estimate_step_0} holds when $\alpha + |\beta| = 0$ by the estimate in \eqref{linear_equation_1_estimates} and the initial data in \eqref{linear_equation_1}. Seeking to apply strong induction, suppose that the statements in \eqref{v_1_estimate_step_0} hold for $0 \leq \alpha + |\beta| \leq m-1$. By differentiating \eqref{linear_equation_1} and rearranging terms, we have that when $\alpha + |\beta| = m$,
\begin{equation} \label{v_1_estimate_step_1_1}
\partial^{\alpha}_{t}\partial^{\beta}_{x}v_{1} = i\Delta\left(\partial^{\alpha-1}_{t}\partial^{\beta}_{x} \right)v_{1} + i\left(\partial^{\alpha-1}_{t}\partial^{\beta}_{x} \right)(pv_{1}) -i \partial^{\alpha-1}_{t}\partial^{\beta}_{x} F.
\end{equation}
Thus, the induction hypothesis implies that $\partial^{\alpha}_{t}\partial^{\beta}_{x}v_{1}(\cdot, 0) = 0$. Again considering the equation obtained by differentiating \eqref{linear_equation_1}, we have (with $\alpha + |\beta| = m$) that
\begin{equation} \label{v_1_estimate_step_1}
(i\partial_{t} + \Delta + p)\left(\partial^{\alpha}_{t}\partial^{\beta}_{x}v_{1} \right) = \\
\partial^{\alpha}_{t}\partial^{\beta}_{x}F + \sum_{\alpha_{0} + |\beta_{0}| = 0}^{m-1} C(\alpha_{0}, \beta_{0})\left(\partial^{\alpha - \alpha_{0}}_{t}\partial^{\beta - \beta_{0}}_{x}p\right)\partial^{\alpha_{0}}_{t}\partial^{\beta_{0}}_{x}v_{1}.
\end{equation}
The claimed estimate in \eqref{v_1_estimate_step_0} now follows from \eqref{linear_equation_1_estimates} and the induction hypothesis. We conclude that 
\begin{equation} \label{v_2_estimate_step_2}
\Vert v_{1} \Vert_{H^{m}(Q)} \leq C \Vert F \Vert_{H^{m}(Q)}
\end{equation}
for a constant $C = C(\Omega, T, \Vert p \Vert_{W^{m,\infty}(Q)}) > 0$, which combined with \eqref{v_2_estimate_step_0} yields \eqref{linear_equation_boundary_data_estimate}. The fact that $\partial^{k}_{t}v(\cdot, 0) = 0$ for integers $0 \leq k \leq m$ is a result of \eqref{linear_equation_1_0} and \eqref{v_1_estimate_step_0}.
\end{proof}
Using Theorem \ref{linear_well-posedness}, we can establish the following result for the well-posedness of \eqref{two_component_GPEqn} when the boundary data are small. Recall Remark \ref{linear_nonlinear_coefficient_assumptions_0}(a) for the reason that we work in the space $H^{3}(Q)$. For $r > 0$, define the neighborhood
\begin{equation} \label{boundary_neighborhood}
\mathcal{B}_{r}(0) \coloneqq \mathcal{H}^{9/2,9/2}(\Sigma) \cap \left\{g \in H^{9/2,9/2}(\Sigma) \,\big|\, \Vert g \Vert_{H^{9/2,9/2}(\Sigma)} \leq r\right\},
\end{equation}
where $\mathcal{H}^{9/2,9/2}(\Sigma)$ was defined in \eqref{boundary_neighborhood_0}.  
\begin{theorem} \label{nonlinear_well-posedness}
Let $\Omega$ be a smooth bounded domain in $\mathbb{R}^{n}\; (n = 2, 3)$ and let $j\in \{1,2\}$. In \eqref{two_component_GPEqn}, suppose that $p_{j} \in W^{3,\infty}(Q)$ and that $q_{j}, \tilde{q}_{j} \in H^{3}(Q)$. Then, there exists $r \in (0,1)$ such that whenever $(g_{1},g_{2}) \in \left(\mathcal{B}_{r}(0)\right)^{2}$, there is a unique solution $(u_{1},u_{2}) \in \left(H^{3}(Q)\right)^{2}$ to \eqref{two_component_GPEqn} satisfying
\begin{equation} \label{nonlinear_equation_boundary_data_estimate}
\Vert u_{j}\Vert_{H^{3}(Q)} \leq C \Vert g_{j}\Vert_{H^{9/2,9/2}(\Sigma)}
\end{equation}
for some $C > 0$. The constants $C$ and $r$ depend on $\Omega, T, \left\Vert p_{j} \right\Vert_{W^{3,\infty}(Q)}, \left\Vert q_{j} \right\Vert_{H^{3}(Q)}, \left\Vert \tilde{q}_{j} \right\Vert_{H^{3}(Q)}$. Moreover, we have that $\partial^{k}_{t}u_{j}(\cdot, 0) = 0$ for $0 \leq k \leq 3$.
\end{theorem}
\begin{proof}
The proof relies on a standard argument based on the Contraction Mapping Principle. We refer the readers to \cite{AZ2025, LLZ2023, LUY2024, LOST2022} for related examples.\\ \\
\textit{Step 1 (A Fixed Point Equation).}\\ \indent
We seek to show that there is a solution of the form $(u_{1}, u_{2}) \coloneqq (\tilde{u}_{1} + w_{1}, \tilde{u}_{2} + w_{2})$ where $\tilde{u}_{j}$ satisfies
\begin{equation} \label{boundary_part_0}
\begin{cases}
\mathcal{T}_{j} \tilde{u}_{j} = 0 \text{ on } Q,\\
\tilde{u}_{j}(x,t) = g_{j}(x,t) \text{ on } \Sigma,\\
\tilde{u}_{j}(x,0) = 0 \text{ on } \Omega \times \{0\}.
\end{cases}
\end{equation}
and $w_{j}$ satisfies
\begin{equation} \label{source_part_0}
\begin{cases}
\mathcal{T}_{j} w_{j}  = f_{j}(u_{1},u_{2}) \coloneqq -q_{j}(x,t)\left|\tilde{u}_{j} + w_{j}\right|^{2}(\tilde{u}_{j} + w_{j})\\
-\; \tilde{q}_{j}(x, t)|\tilde{u}_{\iota(j)} + w_{\iota(j)}|^{2}(\tilde{u}_{j} + w_{j}) \text{ on } Q,\\
w_{j}(x,t) = 0 \text{ on } \Sigma,\\
w_{j}(x,0) = 0 \text{ on } \Omega \times \{0\}.
\end{cases}
\end{equation}
For $g_{j} \in \mathcal{H}^{9/2,9/2}(\Sigma)$ (recall that this space was defined in \eqref{boundary_neighborhood_0}), Theorem \ref{linear_well-posedness} (with $m = 3$) implies that \eqref{boundary_part_0} has a unique solution in $H^{3}(Q)$ satisfying the estimate
\begin{equation} \label{boundary_part_0_Theorem_2.2_estimate}
\Vert \tilde{u}_{j} \Vert_{H^{3}(Q)} \leq C \Vert g_{j} \Vert_{H^{9/2,9/2}(\Sigma)}
\end{equation}
for a constant $C = C\left(\Omega, T, \Vert p_{j} \Vert_{W^{3,\infty}(Q)}\right) > 0$. Moreover, we have that
\begin{equation} \label{boundary_part_0_vanishing_derivatives}
\partial^{k}_{t}\tilde{u}_{j}(\cdot,0) = 0 \quad (0 \leq k \leq 3). 
\end{equation}
Suppose also that
\begin{equation} \label{source_neighborhood_-1}
w_{j} \in \mathcal{H}^{3}(Q) \coloneqq \left\{w \in H^{3}(Q) \,\big|\, \partial^{k}_{t}w(\cdot, 0) = 0 \text{ for } 0 \leq k \leq 3 \right\}.
\end{equation}
With reference to the source term in \eqref{source_part_0}, one then has that
\begin{equation} \label{source_neighborhood_1}
\partial_{t}^{k}f_{j}(u_{1}(\cdot, 0), u_{2}(\cdot, 0)) = 0 \quad (0 \leq k \leq 3)
\end{equation}
and as a result, we can apply the results of Theorem \ref{linear_well-posedness} to \eqref{source_part_0}. In particular, we have the estimate
\begin{equation} \label{source_part_0_Theorem_2.2_estimate}
\Vert w_{j} \Vert_{H^{3}(Q)} \leq C\Vert f_{j} \Vert_{H^{3}(Q)}
\end{equation}
for a constant $C = C\left(\Omega, T, \Vert p_{j} \Vert_{W^{3,\infty}(Q)}\right) > 0$. \\ \indent
For $r > 0$, define the neighborhood
\begin{equation} \label{source_neighborhood}
\mathcal{C}_{r}(0) \coloneqq \mathcal{H}^{3}(Q) \cap \left\{g \in H^{3}(Q) \,\big|\, \Vert g \Vert_{H^{3}(Q)} \leq r\right\}.
\end{equation}
For some $g_{j} \in \mathcal{B}_{r}(0)$ that is fixed and $w_{j} \in \mathcal{C}_{r}(0)$ (where $r > 0$ will be fixed later), consider the unique solution to \eqref{source_part_0}, which we denote by $S_{j}(u_{1}, u_{2}) \coloneqq S(f_{j}(u_{1},u_{2})) \in H^{3}(Q)$. We seek to find a fixed point $(w_{1},w_{2})$ such that
\begin{equation} \label{fixed_point_condition_0}
(w_{1},w_{2}) = S(\tilde{u}_{1} + w_{1}, \tilde{u}_{2} + w_{2}) \coloneqq (S_{1}(\tilde{u}_{1} + w_{1}, \tilde{u}_{2} + w_{2}), S_{2}(\tilde{u}_{1} + w_{1}, \tilde{u}_{2} + w_{2}))
\end{equation}
We will show that \eqref{fixed_point_condition_0} has a fixed point by showing that the solution map
\begin{equation} \label{contraction_map_candidate}
(w_{1},w_{2}) \longmapsto S(\tilde{u}_{1} + w_{1}, \tilde{u}_{2} + w_{2})
\end{equation}
is a contraction map on $(\mathcal{C}_{r}(0))^{2}$ when $r > 0$ is sufficiently small and $g_{j} \in \mathcal{B}_{r}(0)$. \\ \\
\textit{Step 2 (Closure of the Solution Map).}\\ \indent
We begin by checking that if $(g_{1},g_{2}) \in (\mathcal{B}_{r}(0))^{2}$ and $(w_{1}, w_{2})\in (\mathcal{C}_{r}(0))^{2}$ for $r > 0$ sufficiently small, then $S(f_{1},f_{2}) \in (\mathcal{C}_{r}(0))^{2}$. Indeed, as a result of the estimates in \eqref{boundary_part_0_Theorem_2.2_estimate} and \eqref{source_part_0_Theorem_2.2_estimate} along with the fact that $H^{3}(Q)$ is a Banach algebra, we have
\begin{equation} \label{closure_verification}
\Vert S_{j} \Vert_{H^{3}(Q)} \leq C\Vert |\tilde{u}_{j} + w_{j}|^{2}(\tilde{u}_{j} + w_{j}) + |\tilde{u}_{\iota(j)} + w_{\iota(j)}|^{2}(\tilde{u}_{j} + w_{j}) \Vert_{H^{3}(Q)} \leq Cr^{3} \leq r
\end{equation}
for an appropriate choice of $r \in (0, 1)$.\\ \\
\textit{Step 3 (Verifying the Contraction Mapping Property).}\\ \indent
Now, we will show that if $r \in (0, 1)$ is sufficiently small and
\begin{equation} \label{contraction_verification_0_0_0}
\begin{split}
(g_{1},g_{2}) \in (\mathcal{B}_{r}(0))^{2},\\
(w_{1},w_{2}),(W_{1},W_{2}) \in (\mathcal{C}_{r}(0))^{2},
\end{split}
\end{equation}
then there exists $C \in (0, 1)$ such that
\begin{equation} \label{contraction_verification_0_0}
\begin{split}
\Vert S(\tilde{u}_{1} + w_{1}, \tilde{u}_{2} + w_{2}) - S(\tilde{u}_{1} + W_{1}, \tilde{u}_{2} + W_{2}) \Vert_{\left(H^{3}(Q)\right)^{2}}\\
\leq C\left(\left\Vert (w_{1}, w_{2}) - (W_{1}, W_{2}) \right\Vert_{\left(H^{3}(Q)\right)^{2}} \right).
\end{split}
\end{equation}
When we apply Theorem \ref{linear_well-posedness} to \eqref{source_part_0}, we obtain the following estimate for each component of the solution map:
\begin{equation} \label{contraction_verification_0}
\begin{split}
\Vert S_{j}(\tilde{u}_{1} + w_{1}, \tilde{u}_{2} + w_{2}) - S_{j}(\tilde{u}_{1} + W_{1}, \tilde{u}_{2} + W_{2}) \Vert_{H^{3}(Q)}\\
\leq C \Vert f_{j}(\tilde{u}_{1} + w_{1}, \tilde{u}_{2} + w_{2}) - f_{j}(\tilde{u}_{1} + W_{1}, \tilde{u}_{2} + W_{2})  \Vert_{H^{3}(Q)}.
\end{split}
\end{equation}
To further develop the estimate in \eqref{contraction_verification_0}, we recall the definition of $f_{j}$ in \eqref{source_part_0} and note the following for the coupled source term:
\begin{equation} \label{contraction_verification_1}
\begin{split}
\left\Vert |\tilde{u}_{\iota(j)} + w_{\iota(j)}|^{2}(\tilde{u}_{j} + w_{j}) - |\tilde{u}_{\iota(j)} + W_{\iota(j)}|^{2}(\tilde{u}_{j} + W_{j})  \right\Vert_{H^{3}(Q)} \\
\leq \left\Vert |\tilde{u}_{\iota(j)} + w_{\iota(j)}|^{2}(\tilde{u}_{j} + w_{j}) - |\tilde{u}_{\iota(j)} + w_{\iota(j)}|^{2}(\tilde{u}_{j} + W_{j})  \right\Vert_{H^{3}(Q)}\\
+\; \left\Vert |\tilde{u}_{\iota(j)} + w_{\iota(j)}|^{2}(\tilde{u}_{j} + W_{j}) - |\tilde{u}_{\iota(j)} + W_{\iota(j)}|^{2}(\tilde{u}_{j} + W_{j})  \right\Vert_{H^{3}(Q)}\\
\leq Cr^{2} \left\Vert w_{j} - W_{j} \Vert_{H^{3}(Q)} + Cr \Vert |\tilde{u}_{\iota(j)} + w_{\iota(j)}|^{2} - |\tilde{u}_{\iota(j)} + W_{\iota(j)}|^{2} \right\Vert_{H^{3}(Q)}.
\end{split}
\end{equation}
In \eqref{contraction_verification_1}, we used the Banach algebra property of $H^{3}(Q)$ and the assumptions in \eqref{contraction_verification_0_0_0}. Now, by adding and subtracting $(\tilde{u}_{\iota(j)} + w_{\iota(j)})(\overline{\tilde{u}_{\iota(j)} + W_{\iota(j)}})$ inside the last term on the right-hand side of \eqref{contraction_verification_1}, one can show that
\begin{equation} \label{contraction_verification_2}
\left\Vert |\tilde{u}_{\iota(j)} + w_{\iota(j)}|^{2} - |\tilde{u}_{\iota(j)} + W_{\iota(j)}|^{2}  \right\Vert_{H^{3}(Q)} \leq Cr \left\Vert w_{\iota(j)} - W_{\iota(j)} \right\Vert_{H^{3}(Q)}.
\end{equation}
Taking $r \in (0, 1)$ sufficiently small and similarly estimating the terms in \eqref{contraction_verification_0} resulting from the uncoupled source term in \eqref{source_part_0}, the claimed inequality in \eqref{contraction_verification_0_0} follows.\\ \indent
So, the Contraction Mapping Principle implies that there is a unique fixed point $(w_{1}, w_{2}) \in (\mathcal{C}_{r}(0))^{2}$ satisfying \eqref{fixed_point_condition_0} and thus a unique solution of the form $(u_{1}, u_{2}) \coloneqq (\tilde{u}_{1} + w_{1}, \tilde{u}_{2} + w_{2})$ to \eqref{two_component_GPEqn}. From equations \eqref{boundary_part_0_vanishing_derivatives} and \eqref{source_neighborhood_-1} regarding the vanishing derivatives of $\tilde{u}_{j}$ and $w_{j}$, we have the fact that $\partial^{k}_{t}u_{j}(\cdot, 0) = 0$ for $0 \leq k \leq 3$.\\ \\
\textit{Step 4 (Proof of the Stability Estimate).}\\ \indent
From the first inequality in \eqref{closure_verification} and the assumptions that
\begin{equation} \label{stability_estimate_0_1}
\begin{split}
(g_{1},g_{2}) \in (\mathcal{B}_{r}(0))^{2},\\
(w_{1}, w_{2})\in (\mathcal{C}_{r}(0))^{2}
\end{split}
\end{equation}
we can write
\begin{equation}\label{stability_estimate_0}
\Vert w_{j} \Vert_{H^{3}(Q)} = \Vert S_{j} \Vert_{H^{3}(Q)} \leq Cr^{2}\left(\Vert \tilde{u}_{j} \Vert_{H^{3}(Q)} + \Vert w_{j} \Vert_{H^{3}(Q)}\right).
\end{equation}
After taking $r \in (0, 1)$ sufficiently small and applying \eqref{boundary_part_0_Theorem_2.2_estimate} to estimate $\Vert \tilde{u}_{j} \Vert_{H^{3}(Q)}$, equation \eqref{stability_estimate_0} becomes
\begin{equation} \label{stability_estimate_1}
\Vert w_{j} \Vert_{H^{3}(Q)} \leq C\Vert g_{j} \Vert_{H^{9/2,9/2}(\Sigma)}.
\end{equation}
Now, the fact that $u_{j} \coloneqq \tilde{u}_{j} + w_{j}$ yields the stability estimate in \eqref{nonlinear_equation_boundary_data_estimate}. We see from Theorem \ref{linear_well-posedness} and the use of the Banach algebra property to absorb contributions from $\left\Vert q_{j} \right\Vert_{H^{3}(Q)}, \left\Vert \tilde{q}_{j} \right\Vert_{H^{3}(Q)}$ that the constants $C, r$ depend on the parameters $\Omega, T, \left\Vert p_{j} \right\Vert_{W^{3,\infty}(Q)}, \left\Vert q_{j} \right\Vert_{H^{3}(Q)}, \left\Vert \tilde{q}_{j} \right\Vert_{H^{3}(Q)}$.
\end{proof}
\begin{remark} \label{stability_estimate_general_coupled_system}
The estimate \eqref{nonlinear_equation_boundary_data_estimate} is a special case of the more general inequality
\begin{equation} \label{stability_estimate_0_0_alternative_term}
\Vert (u_{1}, u_{2}) \Vert_{\left(H^{3}(Q)\right)^{2}} \leq C\Vert (g_{1}, g_{2}) \Vert_{\left(H^{3}(Q)\right)^{2}}.
\end{equation}
That the stability estimate \eqref{nonlinear_equation_boundary_data_estimate} involves only $g_{j}$ rather than $(g_{1}, g_{2})$ is a consequence of the structure of the nonlinear term $f_{j}$ in \eqref{source_part_0}. In a general coupled system, both $g_{1}$ and $g_{2}$ contribute to the solution. As an example, suppose that we instead considered a system where
\begin{equation} \label{alternative_nonlinear_term}
f_{j}(u_{1}, u_{2}) \coloneqq -q_{j}|u_{j}|^{2}u_{j} + \tilde{q}_{j}u_{\iota(j)}^{3}.
\end{equation}
We will now replace the nonlinearity in \eqref{two_component_GPEqn_small_data} with the one defined by \eqref{alternative_nonlinear_term} and try to prepare an estimate similar to \eqref{stability_estimate_0}. We begin with an estimate similar to \eqref{closure_verification}, namely
\begin{equation} \label{stability_estimate_0_alternative_term}
\begin{split}
\Vert w_{j} \Vert_{H^{3}(Q)} = \Vert S_{j} \Vert_{H^{3}(Q)}
\leq C\Vert |\tilde{u}_{j} + w_{j}|^{2}(\tilde{u}_{j} + w_{j}) + (\tilde{u}_{\iota(j)} + w_{\iota(j)})^{3} \Vert_{H^{3}(Q)}.
\end{split}
\end{equation}
Now, we seek to factor $r^{2}$ as in \eqref{stability_estimate_0} using the assumptions $(g_{1}, g_{2}) \in (\mathcal{B}_{r}(0))^{2}$ and $(w_{1}, w_{2}) \in (\mathcal{C}_{r}(0))^{2}$. This yields
\begin{equation} \label{stability_estimate_1_alternative_term}
\begin{split}
\Vert w_{j} \Vert_{H^{3}(Q)} \leq Cr^{2}\left(\Vert \tilde{u}_{j} \Vert_{H^{3}(Q)} + \Vert w_{j} \Vert_{H^{3}(Q)} + \Vert \tilde{u}_{\iota(j)} \Vert_{H^{3}(Q)} + \Vert w_{\iota(j)} \Vert_{H^{3}(Q)}\right), 
\end{split}
\end{equation}
where we majorized the terms $|\tilde{u}_{j} + w_{j}|^{2}$ and $(\tilde{u}_{\iota(j)} + w_{\iota(j)})^{2}$ from the right-hand side of \eqref{stability_estimate_0_alternative_term} by a constant multiple of $r^{2}$. Applying Theorem \ref{linear_well-posedness} to $\Vert \tilde{u}_{j} \Vert_{H^{3}(Q)}$ and $\Vert \tilde{u}_{\iota(j)} \Vert_{H^{3}(Q)}$, we see from \eqref{stability_estimate_1_alternative_term} that both $g_{1}$ and $g_{2}$ contribute to the solution in this case. Repeating the argument yields the same estimate for $\Vert w_{\iota(j)}\Vert_{H^{3}(Q)}$. So, adding the estimates for $\Vert\tilde{u}_{j} + w_{j}\Vert_{H^{3}(Q)}\;(j = 1,2)$ yields \eqref{stability_estimate_0_0_alternative_term}.
\end{remark}
\section{Linearization for the Coupled Schr\"{o}dinger System} \label{rigorous_linearization}
In this section and the rest of paper, we assume that $\Omega$ is a smooth bounded domain in $\mathbb{R}^{2}$ or $\mathbb{R}^{3}$. Recall that in relation to \eqref{two_component_GPEqn_small_data}, we had taken $j \in \{1,2\}$, $k \in \{1,2\}$, $l \in \{1,2,3\}$, and $\epsilon_{l} \in \mathbb{R}$. \\ \indent
For $f^{(j)}_{l} \in \mathcal{H}^{9/2,9/2}(\Sigma)$, we now define $v^{(l)}_{j}(x,t)$ the unique solution in $H^{3}(Q)$ to
\begin{equation} \label{first_order_linearization_time-dependent}
\begin{cases}
\left(i\partial_t + \Delta + p_{j}(x,t)\right)v^{(l)}_{j}(x,t) = 0 \text{ on } Q,\\
v^{(l)}_{j}(x,t) = f^{(j)}_{l}(x,t) \text{ on } \Sigma,\\
v^{(l)}_{j}(x,0) = 0 \text{ on } \Omega \times \{0\}.
\end{cases}
\end{equation}
By Theorem \ref{linear_well-posedness}, we have that
\begin{equation} \label{first_order_linearization_vanishing_derivatives}
\begin{split}
\left\Vert v^{(l)}_{j} \right\Vert_{H^{3}(Q)} \leq C\left\Vert f^{(j)}_{l} \right\Vert_{H^{9/2,9/2}(\Sigma)},\\
\partial_{t}v^{(l)}_{j}(\cdot,0) = \partial^{2}_{t}v^{(l)}_{j}(\cdot, 0) = 0,
\end{split}
\end{equation}
where $C = C\left(\Omega, T, \Vert p_{j} \Vert_{H^{3}(Q)}\right) > 0$. \\ \indent
Now, given $(\alpha, \beta, \gamma) \in \{1,2,3\}^{3}$, define the functions $\varrho_{(\alpha,\beta,\gamma)}$ and $\pi_{(\alpha,\beta,\gamma)}$:
\begin{equation} \label{kth_order_time-dependent_source_term_positioning}
\begin{split}
\varrho_{(\alpha,\beta,\gamma)}(1) \coloneqq \alpha,\,\varrho_{(\alpha,\beta,\gamma)}(2) \coloneqq \beta,\,\varrho_{(\alpha,\beta,\gamma)}(3) \coloneqq \gamma,\\
\pi_{(\alpha,\beta,\gamma)} \coloneqq \varrho_{(\alpha,\beta,\gamma)} \circ \pi \quad (\pi \in S_{3}).
\end{split}
\end{equation}
We shall then define the source terms
\begin{equation} \label{kth_order_time-dependent_source_terms}
\begin{split}
\mathcal{S}^{(\alpha, \beta, \gamma)}_{j,k}(x,t) \coloneqq -q_{j,k}\sum_{\pi \in S_{3}}v^{\pi_{(\alpha,\beta,\gamma)}(1)}_{j}v^{\pi_{(\alpha,\beta,\gamma)}(2)}_{j}\overline{v^{\pi_{(\alpha,\beta,\gamma)}(3)}_{j}},\\
\tilde{\mathcal{S}}^{(\alpha, \beta, \gamma)}_{j,k}(x,t) \coloneqq -\tilde{q}_{j,k}\sum_{\pi \in S_{3}}v^{\pi_{(\alpha,\beta,\gamma)}(1)}_{j}v^{\pi_{(\alpha,\beta,\gamma)}(2)}_{\iota(j)}\overline{v^{\pi_{(\alpha,\beta,\gamma)}(3)}_{\iota(j)}},
\end{split}
\end{equation}
where $S_{3}$ is the set of all permutations of $\{1,2,3\}$. By \eqref{first_order_linearization_vanishing_derivatives} and Theorem \eqref{linear_well-posedness}, we can define $w^{(\alpha,\beta,\gamma)}_{j,k}(x,t)$ as the unique solution in $H^{3}(Q)$ to
\begin{equation}\label{kth_order_linearization_time-dependent}
\begin{cases}
\left(i\partial_t + \Delta + p_{j}(x,t)\right)w^{(\alpha,\beta,\gamma)}_{j,k}(x,t) = \mathcal{S}^{(\alpha, \beta, \gamma)}_{j,k}(x,t) + \tilde{\mathcal{S}}^{(\alpha, \beta, \gamma)}_{j,k}(x,t)\text{ on } Q,\\
w^{(\alpha,\beta,\gamma)}_{j,k}(x,t) = 0 \text{ on } \Sigma,\\
w^{(\alpha,\beta,\gamma)}_{j,k}(x,0) = 0 \text{ on } \Omega \times \{0\}.
\end{cases}
\end{equation}
Denoting
\begin{equation} \label{solution_series_remainder_equation_-2}
\mathcal{L}_{j}(x,t) \coloneqq \left(\epsilon_{1}v^{(1)}_{j} + \epsilon_{2}v^{(2)}_{j} + \epsilon_{3}v^{(3)}_{j}\right)(x,t),
\end{equation}
we note the following result which is analogous to Proposition $3.3$ from \cite{LLZ2023}.
\begin{lemma} \label{rigorous_linearization_lemma}
Suppose that $p_{j} \in W^{3,\infty}(Q)$ and the nonlinear coefficients in \eqref{two_component_GPEqn_small_data} satisfy \eqref{coefficient_smallness_assumption}. Recalling $\mathcal{B}_{r}(0)$ (defined in \eqref{boundary_neighborhood}) for the $r > 0$ defined by Theorem \ref{nonlinear_well-posedness}, suppose that the nonlinear equation \eqref{two_component_GPEqn_small_data} is equipped with boundary data
\begin{equation}
\epsilon_{1}f^{(j)}_{1} + \epsilon_{2}f^{(j)}_{2} + \epsilon_{3}f^{(j)}_{3} \in \mathcal{B}_{r}(0),
\end{equation}
where $\epsilon_{l} \in \mathbb{R}\;(l = 1,2,3)$ is small enough to ensure that \eqref{two_component_GPEqn_small_data} has a unique solution in $H^{3}(Q)$ satisfying
\begin{equation} \label{two_component_GP_eqn_small_data_vanishing_derivatives}
\partial^{k}_{t}u_{j,k}(\cdot, 0) = 0 \quad (0 \leq k \leq 3)
\end{equation}
and the estimate \eqref{nonlinear_equation_boundary_data_estimate}. Then, the solution to \eqref{two_component_GPEqn_small_data} can be written in the following form:
\begin{equation} \label{nonlinear_solution_Taylor_series_0}
u_{j,k}(x,t) = \mathcal{L}_{j}(x,t) + \frac{1}{6}\left(\sum_{(\alpha,\beta,\gamma) \in \{1,2,3\}^{3}} \epsilon_{\alpha} \epsilon_{\beta} \epsilon_{\gamma} w^{(\alpha,\beta,\gamma)}_{j,k}(x,t)\right) + \mathcal{R}_{j,k}(x,t),
\end{equation}
where $\mathcal{R}_{j,k}$ satisfies the estimate
\begin{equation}\label{Taylor_series_remainder_estimate_0}
\left \Vert \mathcal{R}_{j,k}\right \Vert_{H^{3}(Q)}
\leq C(\left\vert\epsilon_{1}\right\vert + \left\vert\epsilon_{2}\right\vert + \left\vert\epsilon_{3}\right\vert)^{5} \cdot 
\left[\sum_{l = 1}^{3}\left \Vert f^{(j)}_{l} \right \Vert_{H^{9/2,9/2}(\Sigma)}  \right]^{5}
\end{equation}
for some constant $C = C\left(\Omega, T, \left\Vert p_{j} \right\Vert_{W^{3,\infty}(Q)}, M_{j}\right) > 0$.
\end{lemma}
\begin{proof}
From \eqref{first_order_linearization_time-dependent}-\eqref{kth_order_linearization_time-dependent}, we have that if \eqref{nonlinear_solution_Taylor_series_0} holds, then $R_{j,k}$ satisfies
\begin{equation} \label{Taylor_series_remainder_equation}
\begin{cases}
\mathcal{T}_{j}R_{j,k} = -q_{j,k}\left(|u_{j,k}|^{2}u_{j,k} - \left|\mathcal{L}_{j}\right|^{2}\mathcal{L}_{j} \right) -\; \tilde{q}_{j,k}\left(|u_{\iota(j),k}|^{2}u_{j,k} - \left|\mathcal{L}_{\iota(j)}\right|^{2}\mathcal{L}_{j} \right) \text{ on } Q,\\
R_{j,k}(x,t) = 0 \text{ on } \Sigma,\\
R_{j,k}(x,0) = 0 \text{ on } \Omega \times \{0\}.
\end{cases}
\end{equation}
This a result of the fact that
\begin{equation} \label{Taylor_series_remainder_equation_1}
\frac{1}{6}\left(\sum_{(\alpha,\beta,\gamma) \in \{1,2,3\}^{3}} \epsilon_{\alpha} \epsilon_{\beta} \epsilon_{\gamma} \left(\mathcal{S}^{(\alpha,\beta,\gamma)}_{j,k} + \tilde{\mathcal{S}}^{(\alpha,\beta,\gamma)}_{j,k} \right)\right) = -q_{j,k}\left|\mathcal{L}_{j} \right|^{2}\mathcal{L}_{j} - \tilde{q}_{j,k}\left|\mathcal{L}_{\iota(j)} \right|^{2}\mathcal{L}_{j},
\end{equation}
which can be checked by a direct calculation. Since $q_{j,k}, \tilde{q}_{j,k}, u_{j,k}, \mathcal{L}_{j} \in H^{3}(Q)$, the Banach algebra property of $H^{3}(Q)$ implies that the source term in \eqref{Taylor_series_remainder_equation} belongs to $H^{3}(Q)$. In view of \eqref{first_order_linearization_vanishing_derivatives} and \eqref{two_component_GP_eqn_small_data_vanishing_derivatives}, \eqref{Taylor_series_remainder_equation} has a unique solution in $H^{3}(Q)$ by Theorem \ref{linear_well-posedness}. \\ \indent
The calculation leading to \eqref{Taylor_series_remainder_estimate_0} makes use of the fact that $u_{j,k} - \mathcal{L}_{j}$ satisfies \eqref{two_component_GPEqn_small_data} but with trivial initial and boundary data. Then, for example in the first source term in \eqref{Taylor_series_remainder_equation}, we can add and subtract terms to obtain
\begin{equation} \label{solution_series_remainder_equation_-1}
\begin{split}
-q_{j,k}\left(|u_{j,k}|^{2}u_{j,k} - \left|\mathcal{L}_{j}\right|^{2}\mathcal{L}_{j}\right) = -q_{j,k}\Big(|u_{j,k}|^{2}\left(u_{j,k} - \mathcal{L}_{j} \right)\\
+\;  \mathcal{L}_{j}u_{j,k}\left(\overline{u_{j,k}} - \mathcal{L}_{j}\right) + \left(\mathcal{L}_{j}\right)^{2}\left( u_{j,k} - \overline{\mathcal{L}_{j}} \right)\Big).
\end{split}
\end{equation}
The estimates from Theorem \ref{linear_well-posedness}-\ref{nonlinear_well-posedness} say that
\begin{equation} \label{solution_series_remainder_equation_0}
\begin{split}
\left\Vert u_{j,k} \right\Vert_{H^{3}(Q)} \leq C\sum_{l = 1}^{3}\left|\epsilon_{l}\right|\left \Vert f^{(j)}_{l} \right \Vert_{H^{9/2,9/2}(\Sigma)},\\
\left\Vert \mathcal{L}_{j} \right\Vert_{H^{3}(Q)} \leq C\sum_{l = 1}^{3}\left|\epsilon_{l}\right|\left \Vert f^{(j)}_{l} \right \Vert_{H^{9/2,9/2}(\Sigma)},\\
\left\Vert u_{j,k} - \mathcal{L}_{j} \right\Vert_{H^{3}(Q)} \leq C\left(\left\Vert u_{j,k} \right\Vert^{3}_{H^{3}(Q)} + \left\Vert u_{\iota(j),k} \right\Vert^{2}_{H^{3}(Q)}\left\Vert u_{j,k} \right\Vert_{H^{3}(Q)} \right)
\end{split}
\end{equation}
for a constant $C = C\left(\Omega, T, \left\Vert p_{j} \right\Vert_{W^{3,\infty}(Q)}, M_{j}\right) > 0$. The Banach algebra property of $H^{3}(Q)$ was applied to obtain the last inequality in \eqref{solution_series_remainder_equation_0}. We also used \eqref{coefficient_smallness_assumption} to absorb the contribution of the nonlinear coefficient $q_{j,k}$ into the constant. If we use \eqref{solution_series_remainder_equation_0} to estimate the right-hand side of \eqref{solution_series_remainder_equation_-1} and also treat the second source term in \eqref{Taylor_series_remainder_equation} similarly, we obtain the inequality in \eqref{Taylor_series_remainder_estimate_0}.
\end{proof}
As a consequence of \eqref{nonlinear_solution_Taylor_series_0}, we have the following result.
\begin{corollary} \label{nonlinear_solution_Taylor_series_corollary}
Suppose that the conditions of Lemma \ref{rigorous_linearization_lemma} are satisfied. Then, we have that
\begin{equation}
v^{(l)}_{j} = \partial_{\epsilon_{l}}u_{j,k}\big|_{\epsilon_{l} = 0} \quad \text{and} \quad w^{(\alpha,\beta,\gamma)}_{j,k} = \partial^{3}_{\epsilon_{\alpha}\epsilon_{\beta}\epsilon_{\gamma}}u_{j,k}\big|_{\epsilon_{\alpha} = \epsilon_{\beta} = \epsilon_{\gamma} = 0},
\end{equation}
where the partial derivatives are taken in the sense of the mapping from the boundary data $\epsilon_{1}f^{(j)}_{1} + \epsilon_{2}f^{(j)}_{2} + \epsilon_{3}f^{(j)}_{3} \in \mathcal{B}_{r}(0)$ to the unique solution of \eqref{two_component_GPEqn_small_data} in $H^{3}(Q)$.
\end{corollary}
This justifies the definitions of the first and third-order linearizations in \eqref{heuristic_first_order_linearization}-\eqref{heuristic_kth_order_linearization}.
\begin{remark} \label{linear_coefficient_recovery}
In view of Corollary \ref{nonlinear_solution_Taylor_series_corollary}, one can apply Theorem $1.1$ from \cite{BB2020} to \eqref{first_order_linearization_time-dependent} in order to show that the linearized Dirichlet-to-Neumann map defined by
\begin{equation}
f^{(j)}_{l} \in \mathcal{B}_{r}(0) \longmapsto \partial_{\nu}v^{(l)}_{j}\big|_{\Sigma_{j}}
\end{equation}
stably determines the coefficient $p_{j}$.
\end{remark}
\section{Proof of the Inverse Stability Estimate} \label{stability_estimate_proof}
We will now use the third-order linearization \eqref{kth_order_linearization_time-dependent} and the method developed in \cite{LLZ2023} for the recovery of nonlinear coefficients in a single equation to prove the inverse stability estimate in Theorem \ref{main_inverse_stability_theorem}.
\subsection{An Integral Identity} \label{integral_identity_section}
We will now prepare an integral identity that will enable us to derive a Fourier transform estimate for the coefficients in \eqref{two_component_GPEqn_small_data}. In order to prepare the integral identity and to later apply a unique continuation result, we need to restrict the solution of \eqref{kth_order_linearization_time-dependent} to a subset of $\mathcal{O}^{(j)}$ (previously defined in relation to \eqref{near_the_boundary_a_priori_assumption}) and to a subinterval $(0, T^{*})$ where $0 < T^{*} < T$. \\ \indent
To this end, we define the relatively open subsets $\mathcal{O}^{(j)}_{l}\; (l = 1,2,3)$ satisfying
\begin{equation} \label{boundary_neighborhoods_0}
\overline{\mathcal{O}^{(j)}_{l}} \subset \mathcal{O}^{(j)}_{l-1} \subset \mathcal{O}^{(j)}, 
\end{equation}
where we denote $\mathcal{O}^{(j)}_{0} \coloneqq \mathcal{O}^{(j)}$ in the above. We also denote $\Omega_{l} \coloneqq \Omega \setminus \overline{\mathcal{O}^{(j)}_{l}}$. Consider the smooth cut-off function
\begin{equation} \label{smooth_cut-off_function_0}
\chi(x) \coloneqq \begin{cases}1,\: x \in \overline{\Omega} \setminus \mathcal{O}^{(j)}_{2} \\
0,\: x \in \mathcal{O}^{(j)}_{3}
\end{cases}
\end{equation}
satisfying $0 \leq \chi \leq 1$. \\ \indent
Invoking Theorem \ref{linear_well-posedness}, define $v_{0,j}(x,t)$ to be the unique solution in $H^{3}(Q)$ to the adjoint problem
\begin{equation} \label{adjoint_problem_CGO_approach}
\begin{cases}
\left(i\partial_t + \Delta + p_{j}(x,t)\right)v_{0,j} = 0 \text{ on } \Omega \times (0, T^{*}),\\
v_{0,j}(x,t) = f_{0,j}(x,t) \text{ on } \partial \Omega \times (0, T^{*}),\\
v_{0,j}(x,0) = 0 \text{ on } \Omega \times \{T^{*}\},
\end{cases}
\end{equation}
where $f_{0,j} \in \mathcal{H}^{9/2,9/2}(\Sigma)$ will be specified later. 
\begin{proposition} \label{main_integral_identity}
Recalling the functions $\mathcal{S}^{(1,2,3)}_{j,k}, \tilde{\mathcal{S}}^{(1,2,3)}_{j,k}, w^{(1,2,3)}_{j,k} \in H^{3}(Q)$ defined in \eqref{kth_order_linearization_time-dependent}-\eqref{kth_order_time-dependent_source_terms}, consider
\begin{equation} \label{integral_identity_0_0}
\begin{split}
\mathcal{S}^{(1,2,3)}_{j} \coloneqq \mathcal{S}^{(1,2,3)}_{j,1} - \mathcal{S}^{(1,2,3)}_{j,2},\,
\tilde{\mathcal{S}}^{(1,2,3)}_{j} \coloneqq \tilde{\mathcal{S}}^{(1,2,3)}_{j,1} - \tilde{\mathcal{S}}^{(1,2,3)}_{j,2},\\
w^{(1,2,3)}_{j} \coloneqq w^{(1,2,3)}_{j,1} - w^{(1,2,3)}_{j,2}.
\end{split}
\end{equation}
We have that\footnote{Note that $\left[\Delta, \chi \right] = \Delta \chi + 2(\nabla \chi \cdot \nabla)$.}
\begin{equation} \label{integral_identity_0}
\int_{0}^{T^{*}}\int_{\Omega} \left[\chi\left(\mathcal{S}^{(1,2,3)}_{j} + \tilde{\mathcal{S}}^{(1,2,3)}_{j} \right) + \left(\left[\Delta,\chi \right]w^{(1,2,3)}_{j}\right)\right]\overline{v_{0,j}}(x,t) \,dx\,dt = 0.
\end{equation}
\end{proposition}
\begin{proof}
The function $\chi(x)w^{(1,2,3)}_{j}(x,t)$ satisfies \eqref{kth_order_linearization_time-dependent} with source term
\begin{equation} \label{time_dependent_kth-order_cut-off_source_term}
\chi\left(\mathcal{S}^{(1,2,3)}_{j} + \tilde{\mathcal{S}}^{(1,2,3)}_{j}\right)(x,t) + \left(\left[\Delta,\chi \right]w^{(1,2,3)}_{j}\right)(x,t).
\end{equation}
Let us now multiply both sides of the equation satisfied by $\chi w^{(1,2,3)}_{j}$ by $\overline{v_{0,j}}$ and integrate both sides over $\Omega \times (0, T^{*})$. Integrating by parts on the resulting left-hand side, we obtain
\begin{equation} \label{integral_identity_1}
\begin{split}
\int_{\Omega \times (0, T^{*})} \left[(i\partial_{t} + \Delta + p_{j})\left(\chi w^{(1,2,3)}_{j}\right)\right]\overline{v_{0,j}}(x,t)\,dx\,dt \\
 = i\int_{0}^{T^{*}}\int_{\Omega} \partial_{t}\left(\chi w^{(1,2,3)}_{j}\overline{v_{0,j}}\right) \,dx\,dt - i\int_{0}^{T^{*}}\int_{\Omega} \chi w^{(1,2,3)}_{j}\left(\partial_{t}\overline{v_{0,j}}\right)(x,t)\,dx\,dt \\
 +\;\int_{0}^{T^{*}}\int_{\partial \Omega} \partial_{\nu}\left(\chi w^{(1,2,3)}_{j}\right)\overline{v_{0,j}}(x,t)\,dS\,dt - \int_{0}^{T^{*}}\int_{\Omega} \left[\nabla\left(\chi w^{(1,2,3)}_{j}\right) \cdot \nabla \overline{v_{0,j}}\right](x,t)\,dx\,dt\\
 +\;\int_{0}^{T^{*}}\int_{\Omega} \chi w^{(1,2,3)}_{j}\left(p_{j} \overline{v_{0,j}}\right)(x,t)\,dx\,dt.
\end{split}
\end{equation}
The boundary integral in \eqref{integral_identity_1} vanishes since $w^{(1,2,3)}_{j} \equiv 0$ on $\mathcal{O}^{(j)} \times (0, T^{*})$ (owing to the fact that $q_{j,1}-q_{j,2} \equiv \tilde{q}_{j,1} - \tilde{q}_{j,2} \equiv 0$ on $\mathcal{O}^{(j)} \times (0, T^{*})$ from the assumption in \eqref{near_the_boundary_a_priori_assumption}).\\ \indent
Again integrating by parts, note that
\begin{equation} \label{integral_identity_2}
\begin{split}
\int_{0}^{T^{*}} \left[\nabla\left(\chi w^{(1,2,3)}_{j}\right) \cdot \nabla \overline{v_{0,j}}\right](x,t)\,dx\,dt = \int_{0}^{T^{*}}\int_{\partial \Omega} \chi w^{(1,2,3)}_{j}\left(\partial_{\nu}\overline{v_{0,j}} \right)(x,t)\,dS\,dt\\
\;- \int_{0}^{T^{*}}\int_{\Omega} \chi w^{(1,2,3)}_{j}\left(\Delta\overline{v_{0,j}}\right)\,dx\,dt.
\end{split}
\end{equation}
The boundary integral in \eqref{integral_identity_2} vanishes again using the fact that $w^{(1,2,3)}_{j} \equiv 0$ on $\mathcal{O}^{(j)} \times (0, T^{*})$. Then, since $w^{(1,2,3)}_{j}(\cdot,0) \equiv 0$ and $v_{0,j}(\cdot,T^{*}) \equiv 0$, \eqref{integral_identity_1} becomes
\begin{equation} \label{integral_identity_3}
\begin{split}
\int_{\Omega \times (0, T^{*})} \left[(i\partial_{t} + \Delta + p_{j})\left(\chi w^{(1,2,3)}_{j}\right)\right]\overline{v_{0,j}}(x,t)\,dx\,dt\\
=\int_{\Omega \times (0, T^{*})} \left(\chi w^{(1,2,3)}_{j}\right)\left(-i\partial_{t} + \Delta + p_{j} \right)\overline{v_{0,j}}(x,t)\,dx\,dt.
\end{split}
\end{equation}
The right-hand side of \eqref{integral_identity_3} vanishes by \eqref{adjoint_problem_CGO_approach} assuming that $p_{j}$ is real-valued. As a result, we obtain the identity in \eqref{integral_identity_0}. We point out that the calculations in \eqref{integral_identity_1}-\eqref{integral_identity_3} are justified since we had taken $\mathcal{S}^{(1,2,3)}_{j}, \tilde{\mathcal{S}}^{(1,2,3)}_{j}, w^{(1,2,3)}_{j}, v_{0,j} \in H^{3}(Q)$.
\end{proof}
\subsection{CGO Solutions for the Linear Schr\"{o}dinger Equation} \label{CGO_section}
We seek to substitute the CGO solutions for the linear Schr\"{o}dinger equation
\begin{equation} \label{usual_linear_Schrodinger_equation}
(i\partial_{t} + \Delta + p(x,t))u = 0
\end{equation}
(with vanishing data at either $t = 0$ or $t = T^{*}$) constructed in Section $4$ of \cite{LLZ2023} into the integral identity \eqref{integral_identity_0}. Denoting
\begin{equation} \label{CGO_solution_0_0}
\Phi(x,t;\;\omega) \coloneqq \rho(x\cdot \omega - \rho|\omega|^{2}t)
\end{equation}
for a fixed parameter $\rho > 0$ and some fixed $\omega \in \mathbb{R}^{n}\,(n=2,3)$, these solutions are of the form
\begin{equation} \label{CGO_solution_0}
u(x, t) = e^{i\Phi(x,t;\;\omega)}\sum_{k' = 0}^{N} \rho^{-k'}a_{k'}(x,t) + r(x,t),
\end{equation}
where $N > 0$ is a fixed integer. \\ \indent
We introduce a smooth cut-off function
\begin{equation} \label{smooth_cut-off_function_1}
\theta_{h}(t) = \begin{cases}0,\: t \in [0,h] \cup [T^{*} -h, T^{*}]\\
1,\: t \in [2h, T^{*} - 2h]
\end{cases}
\end{equation}
satisfying $0 \leq \theta_{h}(t) \leq 1$, where $h \in (0, \min\{1,T^{*}/4\})$ is chosen to ensure the estimate
\begin{equation} \label{leading_order_amplitude_estimate}
\left\Vert \theta_{h}\right\Vert_{W^{j,\infty}(\mathbb{R})} \leq C_{j}h^{-j}
\end{equation}
for $j \in \mathbb{N}$. Then, in \eqref{CGO_solution_0}, we either take leading order amplitude
\begin{equation} \label{leading_order_amplitude_definition}
a_{0}(x,t) = \theta_{h}(t) \quad \text{ or } \quad a_{0}(x,t) = \theta_{h}(t)e^{i(t\tau + x \cdot \xi)}
\end{equation}
for some fixed $\xi \in \omega^{\perp}$. \\ \indent
In the case of $a_{0} = \theta_{h}(t)$, we take 
\begin{equation} \label{non-adjoint_first_order_amplitude}
a_{1}(x,t) = -(1/2)|\omega|^{-2}(\partial_{t}\theta_{h})(x \cdot \omega).
\end{equation}
Now, taking a fixed $y \in \partial \Omega$, we define
\begin{equation} \label{adjoint_non-adjoint_higher_order_amplitude}
a_{k'}(x + s\omega, t) = \frac{i}{2}\int_{0}^{s}(i\partial_{t} + \Delta + p)a_{k'-1}(x + \tilde{s}\omega, t)\,d\tilde{s}\quad (2 \leq k' \leq N)
\end{equation}
for 
\begin{equation} \label{adjoint_non-adjoint_higher_order_amplitude_1}
x \in L \coloneqq \{x \in \mathbb{R}^{n}\,|\,\omega \cdot (x - y) = 0\}.
\end{equation}
In the case of $a_{0}(x,t) = \theta_{h}(t)e^{i(t\tau + x \cdot \xi)}$, we use \eqref{adjoint_non-adjoint_higher_order_amplitude} to define the higher order amplitudes for $1 \leq k' \leq N$. We observe from \eqref{smooth_cut-off_function_1} that in any case, for any $m \in \mathbb{N}$,
\begin{equation} \label{smooth_cut-off_function_1_observation}
\partial^{\alpha}_{t}a_{k'}(\cdot, 0) = \partial^{\alpha}_{t}a_{k'}(\cdot, T^{*}) = 0\quad(0\leq k' \leq N,\,0 \leq \alpha \leq m).
\end{equation}
The remainder $r(x,t)$ satisfies
\begin{equation} \label{CGO_0_remainder_equation}
\begin{cases}
(i\partial_{t} + \Delta + p)r = -\rho^{-N}e^{i\Phi}(i\partial_{t} + \Delta + p)a_{N} \text{ on } Q,\\
r = 0 \text{ on } \Sigma,\\
r = 0 \text{ on } \Omega \times \{0\},
\end{cases}
\end{equation}
We now establish estimates in $L^{\infty}(Q)$ and $H^{m}(Q)$, where $m \geq 1$ is an integer, for the solution constructed in \eqref{CGO_solution_0}. The $H^{m}(Q)$ estimates that we record have already been proven in Proposition $4.1$ and Remark $4.1$ from \cite{LLZ2023} under the assumption that the linear coefficient $p$ is smooth. We provide the details again to emphasize that it suffices to take $p \in W^{2N + m, \infty}(Q)$, where $N > 0$ is the fixed integer in \eqref{CGO_solution_0}. In the rest of the paper, we shall denote
\begin{equation} \label{Fourier_conjugate_variable_vector}
\left<\tau, \xi \right> = \sqrt{1 + \tau^{2} + |\xi|^{2}}.
\end{equation}
\begin{proposition} \label{LLZ2023_CGO_estimate_proposition}
Fix an integer $N > 0$ in \eqref{CGO_solution_0}, an integer $m \geq 1$, and let $p \in W^{2N + m,\infty}(Q)$. Further consider a fixed vector $\omega \in \mathbb{R}^{n}$ along with fixed parameters $h \in (0, \min\{1,T^{*}/4\})$ and $\tau \in \mathbb{R}$. For $a_{0}(x,t) = \theta_{h}(t)$, we have the estimates
\begin{equation} \label{L_infty_estimates_1}
\begin{split}
\Vert a_{k'} \Vert_{L^{\infty}(Q)} \leq Ch^{-k'}\; (0 \leq k' \leq N),\\
\Vert r \Vert_{L^{\infty}(Q)} \leq C\rho^{-N+2}h^{-N-2}
\end{split}
\end{equation}
and
\begin{equation} \label{H_m_estimates_1}
\begin{split}
\Vert a_{k'} \Vert_{H^{m}(Q)} \leq \tilde{C}h^{-k'-m}\; (0 \leq k' \leq N),\\
\Vert r \Vert_{H^{m}(Q)} \leq \tilde{C}\rho^{-N+2m}h^{-N-m-1}.
\end{split}
\end{equation}
For $a_{0}(x,t) = \theta_{h}(t)e^{i(t\tau + x \cdot \xi)}$, we have the estimates
\begin{equation} \label{L_infty_estimates_0}
\begin{split}
\Vert a_{k'} \Vert_{L^{\infty}(Q)} \leq C\left<\tau, \xi \right>^{2k'}h^{-k'}\; (0 \leq k' \leq N),\\
\Vert r \Vert_{L^{\infty}(Q)} \leq C \rho^{-N+2} \left<\tau,\xi \right>^{2N + 2}h^{-N-2}
\end{split}
\end{equation}
and
\begin{equation} \label{H_m_estimates_0}
\begin{split}
\Vert a_{k'} \Vert_{H^{m}(Q)} \leq \tilde{C} \left<\tau,\xi \right>^{2k'+m}h^{-k'-m}\; (0 \leq k' \leq N),\\
\Vert r \Vert_{H^{m}(Q)} \leq \tilde{C}\rho^{-N+2m} \left<\tau,\xi \right>^{2N + m + 2}h^{-N-m-1}.
\end{split}
\end{equation}
For the $L^{\infty}(Q)$ estimates in \eqref{L_infty_estimates_1} and \eqref{L_infty_estimates_0}, we have $C = C\left(\Omega, T, \left\Vert p \right\Vert_{W^{2N + 1,\infty}(Q)} \right) > 0$. For the $H^{m}(Q)$ estimates in \eqref{H_m_estimates_1} and \eqref{H_m_estimates_0}, we have $\tilde{C} = \tilde{C}(\Omega, T, \left\Vert p \right\Vert_{W^{2N + m,\infty}(Q)}) > 0$.
\end{proposition}
\begin{proof}
Our primary tool will be the estimate of $\Vert \theta_{h} \Vert_{W^{j,\infty}(\mathbb{R})}$ provided in \eqref{leading_order_amplitude_estimate}.\\ \\
\textbf{(I)} We first consider $a_{0} = \theta_{h}(t)$ (whence $a_{1} = -(1/2)|\omega|^{-2}(\partial_{t}\theta_{h})(x \cdot \omega)$ from \eqref{non-adjoint_first_order_amplitude}). In this case, when $k' = 0$ or $1$, the amplitude estimate in \eqref{L_infty_estimates_1} follows directly from \eqref{leading_order_amplitude_estimate}. 
For $2 \leq k' \leq N$, it follows by \eqref{adjoint_non-adjoint_higher_order_amplitude} and induction that the highest order derivative of $p$ arises from $\Delta^{k'-2}p$ while the highest order (nonzero) derivative of $\theta_{h}$ arises from $(i\partial_{t})^{k'-1}a_{1}$. Indeed, if this is true for some $2 \leq k' \leq N$, \eqref{adjoint_non-adjoint_higher_order_amplitude} shows that in $a_{k' + 1}$, the highest order derivative of $p$ is due to $\Delta a_{k'-1}$, thus $\Delta\left(\Delta \right)^{k'-2}p$. Likewise, the highest order (nonzero) derivative of $\theta_{h}$ in $a_{k' + 1}$ is due to $(i\partial_{t})a_{k'-1}$, thus $(i\partial_{t})^{k'}a_{1}$.\\ \indent
Since $h \in (0, \min\{1, T^{*}/4\})$, it follows from \eqref{leading_order_amplitude_estimate} and the preceding derivative count that for $0 \leq k' \leq N$:
\begin{equation} \label{amplitude_estimate_1_step_1}
\begin{split}
\left\Vert a_{k'} \right\Vert_{L^{\infty}(Q)} \leq Ch^{-k'},\\
\left\Vert \partial^{\alpha}_{t}\partial^{\beta}_{x}a_{k'}\right\Vert_{L^{2}(Q)} \leq \tilde{C}h^{-k'-m}
\end{split}
\end{equation}
for $0 \leq \alpha + |\beta| \leq m$ and constants $C = C\left(\Omega, \Vert p \Vert_{W^{2(N-2),\infty}(Q)}\right) > 0, \tilde{C} = \tilde{C}\left(\Omega, \Vert p \Vert_{W^{2(N-2) + m,\infty}(Q)}\right) > 0$. This completes the proofs of the amplitude estimates in \eqref{L_infty_estimates_1}-\eqref{H_m_estimates_1}.\\ \indent
To prove the remainder estimates in \eqref{L_infty_estimates_1}-\eqref{H_m_estimates_1}, we apply Lemma \ref{KS2019_Lemma_2.2_statement} and Theorem \ref{linear_well-posedness} to \eqref{CGO_0_remainder_equation}, which (when combined with the Sobolev Embedding Theorem) yield the inequalities
\begin{equation} \label{remainder_estimate_1_step_0}
\begin{split}
\Vert r \Vert_{L^{\infty}(Q)} \leq C \Vert r \Vert_{C([0, T]; H^{2}(\Omega))} \leq C\rho^{-N}\Vert e^{i\Phi}(i\partial_{t} + \Delta + p)a_{N} \Vert_{H^{1}(0, T; L^{2}(\Omega))},\\
\Vert r \Vert_{H^{m}(Q)} \leq C \rho^{-N} \Vert e^{i\Phi}(i\partial_{t} + \Delta + p)a_{N} \Vert_{H^{m}(Q)}
\end{split}
\end{equation}
for some constant $C = C\left(\Omega, T, \left\Vert p \right\Vert_{W^{m,\infty}(Q)} \right) > 0$. We point that \eqref{smooth_cut-off_function_1_observation} ensured that it was possible to apply Theorem \ref{linear_well-posedness} to \eqref{CGO_0_remainder_equation}. \\ \indent
By substituting \eqref{adjoint_non-adjoint_higher_order_amplitude} (with $k' = N$) into $\Vert e^{i\Phi}(i\partial_{t} + \Delta + p)a_{N} \Vert_{H^{1}(0, T; L^{2}(\Omega))}$, we see that the highest order derivative of $p$ arises from $\partial_{t}(\Delta a_{N})$ (that is, $\partial_{t}\Delta(\Delta^{N - 2}p)$) while the highest order (nonzero) derivative of $\theta_{h}$ arises from $\partial_{t}\left(i\partial_{t}\left(i\partial_{t}\right)^{N-1}a_{1} \right)$. We also obtain a factor of $\rho^{2}$ from $\partial_{t}\left(e^{i\Phi} \right)$. \\ \indent Likewise, in $\Vert e^{i\Phi}(i\partial_{t} + \Delta + p)a_{N} \Vert_{H^{m}(Q)}$, the highest order derivative of $p$ arises from $\partial^{\alpha}_{t}\partial^{\beta}_{x}(\Delta a_{N})\; (\alpha + |\beta| = m)$ (thus from $\partial^{\alpha}_{t}\partial^{\beta}_{x}\Delta(\Delta^{N - 2}p)\; (\alpha + |\beta| = m)$) while the highest order (nonzero) derivative of $\theta_{h}$ arises from $\partial^{m}_{t}\left(i\partial_{t}\left(i\partial_{t}\right)^{N-1}a_{1} \right)$. We obtain a factor of $\rho^{2m}$ from $\partial^{m}_{t}\left(e^{\Phi}\right)$. \\ \indent
These considerations, along with \eqref{leading_order_amplitude_estimate} and the fact that $h \in (0, \min\{1, T^{*}/4\})$, yield the estimates
\begin{equation} \label{remainder_estimate_1_step_1}
\begin{split}
\Vert r \Vert_{L^{\infty}(Q)} \leq C\rho^{-N+2}h^{-(N+2)},\\
\Vert r \Vert_{H^{m}(Q)} \leq \tilde{C}\rho^{-N+2m}h^{-(N + m + 1)}
\end{split}
\end{equation}
for constants $C = C\left(\Omega, T, \Vert p \Vert_{W^{2N - 1, \infty}(Q)}\right) > 0, \tilde{C} = \tilde{C}\left(\Omega, T, \Vert p \Vert_{W^{2(N-1) + m, \infty}(Q)} \right) > 0$. \\ \\
\textbf{(II)} The proofs of \eqref{L_infty_estimates_0}-\eqref{H_m_estimates_0} when $a_{0} = \theta_{h}(t)e^{i(t\tau + x \cdot \xi)}$ are similar except that we obtain factors of $\left<\tau, \xi \right>$ from the exponential term. In what follows, we discuss the proofs of the $\Vert r \Vert_{L^{\infty}(Q)}$ and $\Vert r \Vert_{H^{m}(Q)}$ estimates in \eqref{L_infty_estimates_0}-\eqref{H_m_estimates_0}. \\ \indent 
Let us apply \eqref{adjoint_non-adjoint_higher_order_amplitude} (with $k' = N$ in the case of $a_{0} = \theta_{h}(t)e^{i(t\tau + x \cdot \xi)}$) in the first and second inequalities listed in \eqref{remainder_estimate_1_step_0}. We observe that we obtain factors of $\rho^{2}, \rho^{2m}$ due to $\partial_{t}\left(e^{i\Phi}\right), \partial^{m}_{t}\left(e^{i\Phi}\right)$ in the first and second inequalities from \eqref{remainder_estimate_1_step_0} respectively. In $\Vert e^{i\Phi}(i\partial_{t} + \Delta + p)a_{N} \Vert_{H^{1}(0, T; L^{2}(\Omega))}$, we have that the highest order derivative of $p$ arises from $\partial_{t}(\Delta a_{N})$ (thus from $\partial_{t}\Delta(\Delta^{N-1}p)$ by inductively applying \eqref{adjoint_non-adjoint_higher_order_amplitude}). The highest order derivative of $a_{0}$ with respect to $t$ arises from $\partial_{t}(i\partial_{t} a_{N})$ (thus from $\partial_{t}(\left(i\partial_{t}\right)^{N+1} a_{0})$) while the highest order derivative of $x$ arises from $\Delta a_{N}$ (thus from $\Delta^{N+1} a_{0}$). \\ \indent
In $\Vert e^{i\Phi}(i\partial_{t} + \Delta + p)a_{N} \Vert_{H^{m}(Q)}$, the highest order derivative of $p$ arises from $\partial^{\alpha}_{t}\partial^{\beta}_{x}\Delta a_{N}\; (\alpha + |\beta| = m)$ (thus from $\partial^{\alpha}_{t}\partial^{\beta}_{x}\Delta(\Delta^{N - 1}p)\; (\alpha + |\beta| = m)$). The highest order derivative of $a_{0}$ with respect to $t$ arises from $\partial^{m}_{t}(i\partial_{t} a_{N})$ (thus from $\partial^{m}_{t}(\left(i\partial_{t}\right)^{N+1} a_{0})$) while the highest order derivative with respect to $x$ arises from $\partial^{\beta}_{x}\Delta a_{N}\; (|\beta| = m)$ (thus from $\partial^{\beta}_{x}\Delta^{N+1} a_{0}\; (|\beta| = m)$). \\ \indent
Now, by majorizing each factor of $\tau$ or $\xi$ that appears in these derivatives by $\left< \tau, \xi \right>$ and invoking \eqref{leading_order_amplitude_estimate}, we can develop the estimates in \eqref{remainder_estimate_1_step_0} into
\begin{equation} \label{remainder_estimate_0_step_1}
\begin{split}
\Vert r \Vert_{L^{\infty}(Q)} \leq C \rho^{-N+2} \left<\tau,\xi \right>^{2N + 2}h^{-N-2},\\
\Vert r \Vert_{H^{m}(Q)} \leq \tilde{C}\rho^{-N + 2m}\left<\tau, \xi \right>^{2N + m + 2}h^{-(N + m + 1)}
\end{split}
\end{equation}
for constants $C = C\left(\Omega, T, \Vert p \Vert_{W^{2N + 1,\infty}(Q)} \right) > 0$ and $\tilde{C} = \tilde{C}\left(\Omega, T, \Vert p \Vert_{W^{2N + m,\infty}(Q)} \right)$.
\end{proof}
\begin{remark}
We comment on our motivation for recording $L^{\infty}(Q)$ estimates. We observe from \eqref{H_m_estimates_1} and \eqref{H_m_estimates_0} that an $\mathcal{O}(\rho^{-1})$ estimate for the remainder $r(x,t)$ requires $N \geq 2m + 1$. On the other hand, such a remainder estimate can be fulfilled by taking $N = 3$ in \eqref{L_infty_estimates_1} and \eqref{L_infty_estimates_0}. This smaller choice of $N$ allows us to reduce the regularity required for $p_{j}$ in the forthcoming arguments.
\end{remark}
\subsection{A Fourier Transform Estimate} \label{Fourier_transform_estimate_section}
In what follows for the rest of the paper, we illustrate the proof of the inverse stability estimate \eqref{main_stability_estimate_statement} for $\tilde{q}_{j,k}$, where we fix some $j \in \{1,2\}$. A similar argument holds for the other nonlinear coefficients in \eqref{two_component_GPEqn_small_data}.\\ \indent
We will now substitute CGO solutions of the form \eqref{CGO_solution_0} for $v^{(l)}_{j}$ and $v_{0, j}$ into the integral identity \eqref{integral_identity_0}:
\begin{equation} \label{integral_identity_0_recalled}
\int_{0}^{T^{*}}\int_{\Omega} \left[\chi\left(\mathcal{S}^{(1,2,3)}_{j} + \tilde{\mathcal{S}}^{(1,2,3)}_{j} \right) + \left(\left[\Delta,\chi \right]w^{(1,2,3)}_{j}\right)\right]\overline{v_{0,j}} \,dx\,dt = 0,
\end{equation}
which implies
\begin{equation} \label{integral_identity_0_recalled_1}
\int_{0}^{T^{*}}\int_{\Omega\setminus\mathcal{O}^{(j)}_{3}}\chi\left(\mathcal{S}^{(1,2,3)}_{j} + \tilde{\mathcal{S}}^{(1,2,3)}_{j} \right) \overline{v_{0,j}} \,dx\,dt  \leq \left\vert \int_{0}^{T^{*}}\int_{\Omega} \left(\left[\Delta,\chi \right]w^{(1,2,3)}_{j}\right)\overline{v_{0,j}} \,dx\,dt  \right\vert.
\end{equation}
Recall from \eqref{kth_order_time-dependent_source_terms} and \eqref{integral_identity_0_0} that
\begin{equation} \label{kth_order_time-dependent_source_terms_recalled}
\begin{split}
\mathcal{S}^{(1, 2, 3)}_{j}(x,t) = -2(q_{j,1} - q_{j,2})(x,t)\left[v^{(1)}_{j}v^{(2)}_{j}\overline{v^{(3)}_{j}} + v^{(1)}_{j}v^{(3)}_{j}\overline{v^{(2)}_{j}} + v^{(2)}_{j}v^{(3)}_{j}\overline{v^{(1)}_{j}}\right],\\
\tilde{\mathcal{S}}^{(1, 2, 3)}_{j}(x,t) = -(\tilde{q}_{j,1}-\tilde{q}_{j,2})(x,t)\Big[v^{(1)}_{j}v^{(2)}_{\iota(j)}\overline{v^{(3)}_{\iota(j)}} + v^{(2)}_{j}v^{(1)}_{\iota(j)}\overline{v^{(3)}_{\iota(j)}} + v^{(1)}_{j}v^{(3)}_{\iota(j)}\overline{v^{(2)}_{\iota(j)}} \\
+ v^{(3)}_{j}v^{(1)}_{\iota(j)}\overline{v^{(2)}_{\iota(j)}} + v^{(2)}_{j}v^{(3)}_{\iota(j)}\overline{v^{(1)}_{\iota(j)}} + v^{(3)}_{j}v^{(2)}_{\iota(j)}\overline{v^{(1)}_{\iota(j)}}\Big].
\end{split}
\end{equation}
Also recall the definition of $\Phi(x,t;\;\omega)$ in \eqref{CGO_solution_0_0} and $\theta_{h}(t)$ in \eqref{smooth_cut-off_function_1}. Denoting
\begin{equation}
\begin{split}
a^{(l)}_{j}(x,t) \coloneqq \theta_{h}(t),\\
a_{0,j}(x,t) \coloneqq \theta_{h}(t)e^{i\left(t\tau + x \cdot \xi_{j}\right)},
\end{split}
\end{equation}
we shall fix $N = 3$ and put
\begin{equation} \label{first_order_linearization_CGO_0}
v^{(l)}_{j}(x,t) = e^{i\Phi\left(x,t;\;\omega^{(l)}_{j}\right)}\left(a^{(l)}_{j}(x,t) + \sum_{k' = 1}^{3}\rho^{-k'}a^{(l),j}_{k'}(x,t)\right) + r^{(l)}_{j}(x,t)
\end{equation}
for some choice of $\omega^{(l)}_{j} \in \mathbb{R}^{n}$. For the adjoint IBVP \eqref{adjoint_problem_CGO_approach}, we put
\begin{equation} \label{adjoint_CGO_0}
v_{0,j}(x,t) = e^{i\Phi\left(x,t;\;\omega_{0,j}\right)}\left(a_{0,j}(x,t) +  \sum_{k' = 1}^{3}\rho^{-k'}a^{0,j}_{k'}(x,t)\right) + r_{0,j}(x,t)
\end{equation}
for some choice of $\omega_{0,j} \in \mathbb{R}^{n}$ and $\xi_{j} \in \omega_{0,j}^{\perp}$. With $N = 3$, we observe from Proposition \ref{LLZ2023_CGO_estimate_proposition} that
\begin{equation} \label{amplitude_remainder_estimates_N=3}
\begin{split}
\left \Vert a^{(l)}_{j} \right \Vert_{L^{\infty}(Q)}, \left \Vert a_{0,j} \right \Vert_{L^{\infty}(Q)} \leq C,\\
\left \Vert a^{(l),j}_{k'} \right \Vert_{L^{\infty}(Q)} \leq Ch^{-k'}\; (0 \leq k' \leq 3),\\
\left \Vert r^{(l)}_{j}  \right \Vert_{L^{\infty}(Q)} \leq C\rho^{-1}h^{-5},\\
\left \Vert a^{0,j}_{k'} \right \Vert_{L^{\infty}(Q)} \leq C\left<\tau, \xi_{j} \right>^{2k'}h^{-k'}\; (0 \leq k' \leq 3),\\
\left \Vert r_{0,j} \right \Vert_{L^{\infty}(Q)} \leq C\rho^{-1}\left<\tau, \xi_{j} \right>^{8}h^{-5}
\end{split}
\end{equation}
for a constant $C = C\left(\Omega, T, \left\Vert p_{j} \right\Vert_{W^{7,\infty}(Q)}\right) > 0$.\\ \indent
In what follows, we refer to $a^{(l)}_{j}$ and $a_{0,j}$ as \textit{leading order amplitudes} and the other terms in \eqref{first_order_linearization_CGO_0}-\eqref{adjoint_CGO_0} as \textit{non-leading order amplitudes}. Now, let us stipulate that
\begin{equation} \label{CGO_condition_0}
\omega^{(1)}_{j} + \omega^{(2)}_{\iota(j)} - \omega^{(3)}_{\iota(j)} - \omega_{0,j} = 0 \quad \text{and} \quad
-\left|\omega^{(1)}_{j} \right|^{2} - \left|\omega^{(2)}_{\iota(j)} \right|^{2} + \left|\omega^{(3)}_{\iota(j)} \right|^{2} + \left|\omega_{0,j} \right|^{2} = 0,
\end{equation}
and 
\begin{equation} \label{CGO_condition_1_1_0}
\text{span}\left\{\omega^{(1)}_{j},\omega^{(2)}_{\iota(j)}\right\} \cap \text{span} \left\{\omega^{(2)}_{j},\omega^{(3)}_{j},\omega^{(1)}_{\iota(j)},\omega^{(3)}_{\iota(j)} \right\} = \emptyset.
\end{equation}
To see that this condition can be satisfied, first pick some linearly independent unit vectors $\eta_{1},\eta_{2} \in \mathbb{R}^{n}$. Then, define
\begin{equation} \label{CGO_condition_1_1}
\begin{split}
\omega^{(1)}_{j} \coloneqq \eta_{1},\, \omega^{(2)}_{\iota(j)} \coloneqq -\eta_{1},\\
\omega^{(2)}_{j} = \omega^{(3)}_{j} = \omega^{(1)}_{\iota(j)} = \omega^{(3)}_{\iota(j)} \coloneqq \eta_{2},\, \omega_{0,j} = -\eta_{2}.
\end{split}
\end{equation}
This simultaneously ensures that the exponential phase for the $v^{(1)}_{j}v^{(2)}_{\iota(j)}\overline{v^{(3)}_{\iota(j)}}\overline{v_{0,j}}$ term in left-hand side of \eqref{integral_identity_0_recalled_1} vanishes while the exponential phases for the other leading order amplitude products (that result from substituting \eqref{first_order_linearization_CGO_0}-\eqref{adjoint_CGO_0} into the left-hand side of \eqref{integral_identity_0_recalled_1}) are nonzero.\\ \indent
Denoting $\Phi\left(\omega^{(l)}_{j}\right) \coloneqq \Phi\left(x,t;\;\omega^{(l)}_{j}\right)$, we have for the $v^{(1)}_{j}v^{(2)}_{\iota(j)}\overline{v^{(3)}_{\iota(j)}}\overline{v_{0,j}}$ term that
\begin{equation} \label{CGO_condition_0_0}
\Phi\left(\omega^{(1)}_{j}\right) + \Phi\left(\omega^{(2)}_{\iota(j)}\right) - \Phi\left(\omega^{(3)}_{\iota(j)}\right) - \Phi\left(\omega_{0,j}\right) = 0.
\end{equation}
This results in the leading order amplitudes for this term combining to yield
\begin{equation} \label{leading_order_amplitude_integral_identity_0}
\int_{0}^{T^{*}}\int_{\Omega\setminus\mathcal{O}^{(j)}_{3}} \chi\left(\tilde{q}_{j,2} - \tilde{q}_{j,1} \right)\theta_{h}^{4} e^{-i(\tau t + x \cdot \xi_{j})}\,dx\,dt.
\end{equation}
Now, we seek to isolate the above term and show that the other terms in the left-hand side of \eqref{integral_identity_0_recalled_1} are of order $\rho^{-1}$. As an example of how such terms are of order $\rho^{-1}$, let us consider the product of the leading order amplitudes in the $v^{(2)}_{j}v^{(1)}_{\iota(j)}\overline{v^{(3)}_{\iota(j)}}\overline{v_{0,j}}$ term from \eqref{integral_identity_0_recalled_1}. Denoting
\begin{equation} \label{another_leading_order_amplitude_0_0}
\omega' \coloneqq \omega^{(2)}_{j} + \omega^{(1)}_{\iota(j)} - \omega^{(3)}_{\iota(j)} - \omega_{0,j} \quad \text{and} \quad W \coloneqq \left|\omega^{(2)}_{j} \right|^{2} + \left|\omega^{(1)}_{\iota(j)} \right|^{2} - \left|\omega^{(3)}_{\iota(j)} \right|^{2} - \left|\omega_{0,j} \right|^{2},
\end{equation}
we obtain
\begin{equation} \label{another_leading_order_amplitude_0}
\int_{0}^{T^{*}}\int_{\Omega\setminus\mathcal{O}^{(j)}_{3}} \chi\left(\tilde{q}_{j,2} - \tilde{q}_{j,1} \right)\theta_{h}^{4} e^{i\rho\left(x \cdot \omega' - \rho Wt\right)}e^{-i(\tau t + x \cdot \xi_{j})} \,dx\,dt
\end{equation}
for the product of the leading order amplitudes in the $v^{(2)}_{j}v^{(1)}_{\iota(j)}\overline{v^{(3)}_{\iota(j)}}\overline{v_{0,j}}$ term. Observe that we have $\omega' \neq 0$ because of \eqref{CGO_condition_1_1_0}-\eqref{CGO_condition_1_1}. \\ \indent
Now, denoting $\nabla' = (\partial_{x_{1}},\ldots,\partial_{x_{n}},\partial_{t})$, we observe (as noted in Section $5$ of \cite{LUY2024}) that
\begin{equation} \label{observation_0}
e^{i\rho\left(x \cdot \omega' - \rho Wt \right)} = \nabla'e^{i\rho\left(x \cdot \omega' - \rho W t \right)} \cdot \frac{\left(-i\rho \omega',i\rho^{2}W \right)}{\left|\left(\rho \omega',\rho^{2} W \right)\right|^{2}}.
\end{equation}
Substituting \eqref{observation_0} into \eqref{another_leading_order_amplitude_0}, we may integrate by parts to obtain
\begin{equation} \label{another_leading_order_amplitude_1}
\begin{split}
\int_{\partial((\Omega\setminus\mathcal{O}^{(j)}_{3}) \times (0,T^{*}))} \left(\nu(x,t) \cdot \frac{\left(-i\rho \omega',i\rho^{2} W \right)}{\left|\left(\rho \omega',\rho^{2} W \right)\right|^{2}} \right)\chi\left(\tilde{q}_{j,2} - \tilde{q}_{j,1} \right) \left[\theta_{h}(t) \right]^{4} e^{i\rho\left(x \cdot \omega' - \rho W t\right)}e^{-i(\tau t + x \cdot \xi_{j})}\,dS\\
-\int_{(\Omega\setminus\mathcal{O}^{(j)}_{3}) \times (0, T^{*})}\frac{\left(-i\rho \omega',i\rho^{2} W \right)}{\left|\left(\rho \omega',\rho^{2} W \right)\right|^{2}} \cdot \nabla'\left[\chi\left(\tilde{q}_{j,2} - \tilde{q}_{j,1} \right) \left[\theta_{h}(t) \right]^{4}e^{-i(\tau t + x \cdot \xi_{j})} \right]e^{i\rho\left(x \cdot \omega' - \rho W t\right)}\,dx\,dt.
\end{split}
\end{equation}
Due to the assumption that $\tilde{q}_{j,1} \equiv \tilde{q}_{j,2}$ on $\mathcal{O}^{(j)} \times (0, T)$ and the definition of $\theta_{h}(t)$ in \eqref{smooth_cut-off_function_1}, the boundary integral in \eqref{another_leading_order_amplitude_1} vanishes. Denoting the non-boundary term in \eqref{another_leading_order_amplitude_1} by $I$, we have from the fact that $\omega' \neq 0$ that
\begin{equation} \label{another_leading_order_amplitude_2}
|I| \leq \frac{\rho^{-1}}{\left|(\omega', 0) \right|}\int_{(\Omega\setminus\mathcal{O}^{(j)}_{3}) \times (0, T^{*})} \left| \nabla'\left[\chi\left(\tilde{q}_{j,2} - \tilde{q}_{j,1} \right) \left[\theta_{h}(t) \right]^{4}e^{-i(\tau t + x \cdot \xi_{j})} \right] \right|\,dx\,dt
\end{equation}
By \eqref{coefficient_smallness_assumption} and \eqref{leading_order_amplitude_estimate}, we can further say that
\begin{equation} \label{another_leading_order_amplitude_bound_0}
|I| \leq C\rho^{-1}\left(h^{-1} + \left<\tau, \xi_{j} \right>\right),
\end{equation}
which is what we sought to show. In \eqref{another_leading_order_amplitude_bound_0}, the factor of $h^{-1}$ resulted from applying \eqref{leading_order_amplitude_estimate} to estimate $|\nabla'\theta_{h}|$, the factor $\left<\tau, \xi_{j} \right>$ came from $\nabla'e^{-i(\tau t + x \cdot \xi_{j})}$, and the contribution from $\nabla'\left(\tilde{q}_{j,2} - \tilde{q}_{j,1}\right)$ was absorbed into the constant by \eqref{coefficient_smallness_assumption}. We also absorbed the contribution of $\nabla'\chi$, which depends on $\mathcal{O}^{(j)}_{2}$ and $\mathcal{O}^{(j)}_{3}$ in view of \eqref{smooth_cut-off_function_0}, into the constant. This leads to the parameter dependence $C = C(\omega', \mathcal{O}^{(j)}, T^{*}, M_{j}) > 0$ in \eqref{another_leading_order_amplitude_bound_0}.
\begin{remark} \label{vectors_for_other_coefficients}
Henceforth, we will indicate the dependence of the constant $C$ on vector combinations resulting from the product of exponential phases (for example the expression $|(\omega', 0)|$ in \eqref{another_leading_order_amplitude_2}) by writing $\tilde{C}_{j}$, with the understanding that we have chosen a specific set of vectors for the purpose of recovering the coefficient $\tilde{q}_{j}$. If we were concerned with recovering the coefficient $q_{j}$ for a fixed $j \in \{1,2\}$, one could stipulate
\begin{equation}
\omega^{(1)}_{j} + \omega^{(2)}_{j} - \omega^{(3)}_{j} - \omega_{0,j} = 0 \quad \text{and} \quad
-\left|\omega^{(1)}_{j} \right|^{2} - \left|\omega^{(2)}_{j} \right|^{2} + \left|\omega^{(3)}_{j} \right|^{2} + \left|\omega_{0,j} \right|^{2} = 0,
\end{equation}
which ensures that the exponential phase in the $v^{(1)}_{j}v^{(2)}_{j}\overline{v^{(3)}_{j}}\overline{v_{0,j}}$ term vanishes. In order to ensure that the exponential phase for the other terms doesn't vanish, we can additionally stipulate a condition similar to \eqref{CGO_condition_1_1_0} and then pick
\begin{equation}
\begin{split}
\omega^{(1)}_{j} \coloneqq \eta_{1},\, \omega^{(2)}_{j} \coloneqq -\eta_{1},\\
\omega^{(2)}_{\iota(j)} = \omega^{(3)}_{j} = \omega^{(1)}_{\iota(j)} = \omega^{(3)}_{\iota(j)} \coloneqq \eta_{2},\, \omega_{0,j} = -\eta_{2}
\end{split}
\end{equation}
for some linearly independent vectors $\eta_{1},\eta_{2} \in \mathbb{R}^{n}$.
\end{remark}
Next, we seek to estimate terms that result from products with non-leading order amplitudes after substituting \eqref{first_order_linearization_CGO_0}-\eqref{adjoint_CGO_0} into \eqref{integral_identity_0_recalled_1}. One may use the estimates in \eqref{amplitude_remainder_estimates_N=3} to majorize such terms by expressions of the form
\begin{equation} \label{non-leading_order_amplitude_product_0}
C\rho^{-1}h^{-\alpha}\left<\tau, \xi_{j} \right>^{\beta}
\end{equation}
for $\alpha,\beta > 0$ sufficiently large. Indeed, we see from the exponents in \eqref{amplitude_remainder_estimates_N=3} that all such products can be majorized in magnitude by
\begin{equation} \label{non-leading_order_amplitude_product_1}
C\left(\rho^{-1}h^{-5}\left<\tau, \xi_{j} \right>^{8}\right)^{4}
\end{equation}
for a constant $C = C\left(\Omega, T, T^{*}, \left\Vert p_{j} \right\Vert_{W^{7,\infty}(Q)}, M_{j}\right) > 0$. Combining the considerations in \eqref{another_leading_order_amplitude_0_0}-\eqref{CGO_condition_1_1}, Remark \ref{vectors_for_other_coefficients}, and \eqref{non-leading_order_amplitude_product_0}-\eqref{non-leading_order_amplitude_product_1}, we can develop \eqref{integral_identity_0_recalled_1} into the estimate
\begin{equation} \label{initial_Fourier_estimate_0}
\begin{split}
\int_{0}^{T^{*}}\int_{\Omega\setminus\mathcal{O}^{(j)}_{3}} \chi\left(\tilde{q}_{j,2} - \tilde{q}_{j,1} \right)\theta_{h}^{4} e^{-i(\tau t + x \cdot \xi_{j})}\,dx\,dt\\
\leq \tilde{C}_{j}\rho^{-1}h^{-20}\left<\tau, \xi_{j} \right>^{32} + \left\vert \int_{0}^{T^{*}}\int_{\Omega} \left(\left[\Delta,\chi \right]w^{(1,2,3)}_{j}\right)\overline{v_{0,j}} \,dx\,dt  \right\vert
\end{split}
\end{equation}
for some constant $\tilde{C}_{j} = \tilde{C}_{j}(\Omega, \mathcal{O}^{(j)}, T, T^{*}, \left\Vert p_{j} \right\Vert_{W^{7,\infty}(Q)}, M_{j}) > 0$. Extending $\chi, \tilde{q}_{j,k}$, and $\theta_{h}$ to be zero outside $Q$, we obtain the estimate
\begin{equation} \label{initial_Fourier_estimate_0_0}
\mathcal{F}\left[\chi\left(\tilde{q}_{j,2} - \tilde{q}_{j,1} \right)\theta_{h}^{4}\right] \leq \tilde{C}_{j}\rho^{-1}h^{-20}\left<\tau, \xi_{j} \right>^{32} + \left\vert \int_{0}^{T^{*}}\int_{\Omega} \left(\left[\Delta,\chi \right]w^{(1,2,3)}_{j}\right)\overline{v_{0,j}} \,dx\,dt  \right\vert,
\end{equation}
where ``$\mathcal{F}$" denotes the Fourier transform in the conjugate variables $(\tau, \xi_{j})$. \\ \indent
We now seek to introduce the Neumann data through the commutator term $\left(\left[\Delta,\chi \right]w^{(1,2,3)}_{j}\right)$ appearing in \eqref{initial_Fourier_estimate_0_0}. For this, we will use the following unique continuation result. To state the result, we need to recall the subsets $\mathcal{O}^{(j)}_{l} \subset \mathcal{O}^{(j)}\; (l = 1, 2, 3)$ (and $\Omega_{l} = \Omega \setminus \overline{\mathcal{O}^{(j)}_{l}}$) defined in \eqref{boundary_neighborhoods_0}.
\begin{proposition} \label{unique_continuation_result_BC2009}
Suppose $\mathcal{O} \subset \Omega$ is a relatively open neighborhood of $\partial\Omega$ such that $\partial\mathcal{O} \setminus \partial\Omega$ is $C^{2}$. Let $p(x,t) \in L^{\infty}(Q)$ with $p(x,t) = p(x)$ on $\mathcal{O} \times (0, T)$ and $\Gamma$ be a fixed relatively open subset of $\partial \Omega$. Suppose $u(x,t) \in H^{1}(0, T; H^{1}(\Omega))$ is a solution to
\begin{equation} \label{kth_order_linearization_LLZ2023}
\begin{cases}
(i\partial_{t} + \Delta + p(x,t))u = g(x,t) \text{ on } Q,\\
u(x,t) = 0 \text{ on } \Sigma,\\
u(x,0) = 0 \text{ on } \Omega \times \{0\},
\end{cases}
\end{equation}
where $g(x,t) \in L^{2}(Q)$ is supported in $(\Omega \setminus \mathcal{O}) \times (0, T)$. Then, for any $T^{*} \in (0, T)$, there exist constants $C = C(\Omega, \Gamma, \mathcal{O}, T^{*}, T) > 0, \gamma_{*} > 0$, and $m > 0$ such that
\begin{equation} \label{unique_continuation_statement_BC2009}
\left \Vert u \right \Vert_{L^{2}\left(0, T^{*};H^{1}\left(\Omega_{3}\setminus\Omega_{2}\right)\right)} \leq C \left(\gamma^{-1} + e^{m\gamma} \left \Vert \partial_{\nu}u \right \Vert_{L^{2}(\Gamma \times (0, T))} \right)
\end{equation}
for any $\gamma \geq \gamma_{*}$.
\end{proposition}
\begin{remark} \label{BC2009_unique_continuation_modifications}
\hfill
\begin{itemize}
\item[(a)] Proposition \ref{unique_continuation_result_BC2009} can be proved by following the proof of Lemma $4.6$ in \cite{BC2009}. The proof is based on considering the linear Schr\"{o}dinger equation satisfied by $\chi u$ on $\mathcal{O} \times (0, T)$, where $\chi$ is defined by \eqref{smooth_cut-off_function_0}. The proof also relies on a Carleman estimate with weight function defined on the interval $(-T, T) (\ni t)$. So, we need to extend $\chi u$ to the interval $(-T, T)$, which can be done by defining
\begin{equation*}
u(x,t) = \overline{u(x,-t)} \quad \text{and} \quad g(x,t) = \overline{g(x,-t)}
\end{equation*}
for $(x,t) \in \mathcal{O} \times (-T, 0)$.  Since $p(x,t) = p(x)$ and $g(x,t) = 0$ on $\mathcal{O} \times (0, T)$, the equation satisfied by $\chi u$ then reduces to equation $(4.13)$ in \cite{BC2009}. The rest of the proof then follows as in equations $(4.12)-(4.51)$ in \cite{BC2009}. In fact, the proof in \cite{BC2009} holds under the weaker condition that $\mathcal{O}$ is a relatively open neighborhood of $\Gamma$ rather than the entire boundary. However, the fact that $\mathcal{O}$ is a neighborhood of the entire boundary is essential for other steps towards the proof of the main result in \cite{BC2009} and our proof of Theorem \ref{main_inverse_stability_theorem}. In particular, we need to know that $q_{j,k}$ and $\tilde{q}_{j,k}$ are determined near the entire boundary in order to derive the integral identity \eqref{integral_identity_0}. The reason for the regularity requirement of $\partial\mathcal{O} \setminus \partial\Omega$ is explained in Remark $4.5$ in \cite{AZ2025}.
\item[(b)] Strictly speaking, the inequality in Lemma $4.6$ of \cite{BC2009} estimates $\left \Vert u \right \Vert_{L^{2}\left(0, T_{1};H^{1}\left(\Omega_{3}\setminus\Omega_{2}\right)\right)}$ for a specific choice of $T_{1} \in (0, T)$ such that the constant in \eqref{unique_continuation_statement_BC2009} is independent of $T_{1}$. However, the proof in \cite{BC2009} still holds for any choice of $T^{*} \in (0, T)$ and in this case, the constant that results in \eqref{unique_continuation_statement_BC2009} now depends on $T^{*}$. Also, the estimate in \cite{BC2009} is in terms of $\left \Vert \partial_{\nu}u \right \Vert_{L^{2}(\Gamma \times (-T, T))}$ on the right-hand side. However, this can be majorized by a constant times $\left \Vert \partial_{\nu}u \right \Vert_{L^{2}(\Gamma \times (0, T))}$ in view of how we defined the extended solution in (a).
\item[(c)] Regarding \eqref{unique_continuation_statement_BC2009}, the dependence of the constant on $\Gamma$ is not explicitly mentioned in Lemma $2.1$ from \cite{BC2009}. However, such a dependence is apparent from the Carleman estimate (Lemma $4.2$ in \cite{BC2009}) used to establish the Lemma $2.1$. We refer the reader to Remark $4.6$ in \cite{AZ2025} for further details.
\end{itemize}
\end{remark}
\begin{remark} \label{real_valued_data_assumptions}
We seek to apply the estimate in \eqref{unique_continuation_statement_BC2009} to the commutator term in \eqref{initial_Fourier_estimate_0_0}. In view of the discussion in Remark \ref{BC2009_unique_continuation_modifications}(a), this requires extending $w^{(1,2,3)}_{j}$ to $\mathcal{O}^{(j)} \times (-T, T)$. We can use the extension defined in Remark \ref{BC2009_unique_continuation_modifications}(a) if $p_{j}, q_{j,k}, \tilde{q}_{j,k}$ are real-valued and if $v^{(l)}_{j}$ (which defines the source terms in \eqref{kth_order_time-dependent_source_terms}) satisfies
\begin{equation*}
v^{(l)}_{j}(x,t) = \overline{v^{(l)}_{j}(x,-t)}
\end{equation*}
for $(x,t) \in \mathcal{O}^{(j)} \times (-T, 0)$. In turn, this requires the boundary data $f^{(j)}_{l}$ to be real-valued.
\end{remark}
Let us now estimate the commutator term in \eqref{initial_Fourier_estimate_0}:
\begin{equation} \label{initial_Fourier_estimate_1}
\begin{split}
\left|\int_{0}^{T^{*}}\int_{\Omega}\left(\left[\Delta,\chi \right]w^{(1,2,3)}_{j}\right)\overline{v_{0,j}}\,dx\,dt \right| \\
\leq C\left \Vert \left[\Delta,\chi\right]w^{(1,2,3)}_{j} \right\Vert_{L^{2}(\Omega \times (0, T^{*}))} \left \Vert v_{0,j} \right \Vert_{L^{\infty}(\Omega \times (0, T^{*}))} \\
\leq C\left \Vert w^{(1,2,3)}_{j} \right\Vert_{L^{2}(0, T^{*}; H^{1}(\Omega_{3}\setminus\Omega_{2}))} \left \Vert v_{0,j} \right \Vert_{L^{\infty}(\Omega \times (0, T^{*}))} \\
\leq C\left(\gamma_{j}^{-1} + e^{m_{j}\gamma_{j}} \left \Vert \partial_{\nu}w^{(1,2,3)}_{j} \right \Vert_{L^{2}(\Sigma_{j})} \right)\left \Vert v_{0,j} \right \Vert_{L^{\infty}(\Omega \times (0, T^{*}))}\;(\gamma_{j} \geq \gamma_{*,j})
\end{split}
\end{equation}
for some constants $C = C(\Omega, \Gamma_{j}, \mathcal{O}^{(j)}, T^{*}, T) > 0, \gamma_{*,j} > 0, m_{j} > 0$. The last line resulted from substituting the unique continuation estimate \eqref{unique_continuation_statement_BC2009} from Proposition \ref{unique_continuation_result_BC2009}.\\ \indent
Next, observe from \eqref{adjoint_CGO_0} and \eqref{amplitude_remainder_estimates_N=3} that
\begin{equation} \label{initial_Fourier_estimate_1_0}
\begin{split}
\left \Vert v_{0,j} \right \Vert_{L^{\infty}(\Omega \times (0, T^{*}))} \leq \left\Vert a_{0,j} \right\Vert_{L^{\infty}(\Omega \times (0, T^{*}))} + \sum_{k' = 1}^{3}\rho^{-k'}\left \Vert a^{0,j}_{k'} \right \Vert_{L^{\infty}(\Omega \times (0, T^{*}))} + \left \Vert r_{0,j} \right \Vert_{L^{\infty}(\Omega \times (0, T^{*}))}\\
\leq C\left(1 + \rho^{-1}h^{-5}\left<\tau,\xi_{j} \right>^{8} \right),
\end{split}
\end{equation}
where $C = C\left(\Omega, T, \Vert p \Vert_{W^{7,\infty}(Q)}\right) > 0$. In \eqref{initial_Fourier_estimate_1_0}, we majorized each term by $\rho^{-1}h^{-\alpha}\left<\tau,\xi_{j} \right>^{\beta}$ for $\alpha,\beta > 0$ sufficiently large. Substituting \eqref{initial_Fourier_estimate_1_0} into \eqref{initial_Fourier_estimate_1} and then substituting \eqref{initial_Fourier_estimate_1} into \eqref{initial_Fourier_estimate_0_0} yields the following Fourier transform estimate:
\begin{equation} \label{initial_Fourier_estimate_2}
\begin{split}
\mathcal{F}\left[\chi\left(\tilde{q}_{j,2} - \tilde{q}_{j,1} \right)\theta_{h}^{4}\right] 
\leq \tilde{C}_{j}\rho^{-1}h^{-20}\left<\tau, \xi_{j} \right>^{32}\\
+\; \tilde{C}_{j}\left(\gamma_{j}^{-1} + e^{m_{j}\gamma_{j}} \left \Vert \partial_{\nu}w^{(1,2,3)}_{j} \right \Vert_{L^{2}(\Gamma \times (0, T))} \right)\left(1 + \rho^{-1}h^{-5}\left<\tau,\xi_{j} \right>^{8} \right)
\end{split}
\end{equation}
for some constant $\tilde{C}_{j} = \tilde{C}_{j}\left(\Omega, \Gamma_{j}, \mathcal{O}^{(j)}, T, T^{*}, \left\Vert p_{j} \right\Vert_{W^{7,\infty}(Q)} , M_{j}\right) > 0$.
\subsection{Completing the Proof} \label{completing_the_proof_section}
In \eqref{initial_Fourier_estimate_2}, let us fix $\gamma_{j} > 1$ and choose $\rho \coloneqq \gamma_{j}$. Also, pick $\tilde{m}_{j} \geq m_{j} > 0$ such that
\begin{equation} \label{parameterization_2}
e^{m_{j}\gamma_{j}}\rho^{-1} = e^{m_{j}\gamma_{j}}\gamma_{j}^{-1} \leq e^{\tilde{m}_{j}\gamma_{j}}.
\end{equation}
With these choices, \eqref{initial_Fourier_estimate_2} becomes
\begin{equation} \label{initial_Fourier_estimate_3}
\begin{split}
\mathcal{F}\left[\chi\left(\tilde{q}_{j,2} - \tilde{q}_{j,1} \right)\theta_{h}^{4}\right] 
\leq \tilde{C}_{j}\gamma_{j}^{-1}h^{-20}\left<\tau, \xi_{j} \right>^{32}\\
+\; \tilde{C}_{j}\left(\gamma_{j}^{-1} + e^{\tilde{m}_{j}\gamma_{j}} \left \Vert \partial_{\nu}w^{(1,2,3)}_{j} \right \Vert_{L^{2}(\Sigma_{j})} \right)\left(1 + \gamma_{j}^{-1}h^{-5}\left<\tau,\xi_{j} \right>^{8} \right)\\
\leq \tilde{C}_{j}\left(\gamma_{j}^{-1} + e^{\tilde{m}_{j}\gamma_{j}} \left \Vert \partial_{\nu}w^{(1,2,3)}_{j} \right \Vert_{L^{2}(\Sigma_{j})}\right)h^{-20}\left<\tau, \xi_{j} \right>^{32}.
\end{split}
\end{equation}
We will now use \eqref{initial_Fourier_estimate_3} to estimate the $H^{-1}(\mathbb{R}^{n + 1})$-norm and subsequently use interpolation to pass to the $L^{2}(\Omega \times (0, T^{*}))$-norm. To begin with, we have the following for some constant $M > 1$ that will be fixed later:
\begin{equation} \label{negative_one_norm_step_1_0}
\begin{split}
\left \Vert \chi \left(\tilde{q}_{j,1} - \tilde{q}_{j,2} \right) \theta^{4}_{h} \right \Vert^{2}_{H^{-1}(\mathbb{R}^{n + 1})} = \int_{\left<\tau, \xi_{j} \right> \leq M} \left<\tau, \xi_{j} \right>^{-2} \left|\mathcal{F}\left[\chi\left(\tilde{q}_{j,1} - \tilde{q}_{j,2} \right) \theta^{4}_{h}\right]\right|^{2}\,d\tau\,d\xi_{j} \\
+\; \int_{\left<\tau, \xi_{j} \right> > M} \left<\tau, \xi_{j} \right>^{-2} \left|\mathcal{F}\left[\chi\left(\tilde{q}_{j,1} - \tilde{q}_{j,2} \right) \theta^{4}_{h}\right]\right|^{2}\,d\tau\,d\xi_{j}.
\end{split}
\end{equation}
Denote the first and second integrals in \eqref{negative_one_norm_step_1} by $I_{1}$ and $I_{2}$ respectively. Upon substituting \eqref{initial_Fourier_estimate_3} into $I_{1}$, we have that
\begin{equation} \label{negative_one_norm_step_2}
\begin{split}
\vert I_{1} \vert \leq CM^{n + 1} \cdot \tilde{C}_{j}\left(\gamma_{j}^{-1} +  e^{\tilde{m}_{j}\gamma_{j}}\left\Vert \partial_{\nu}w^{(1,2,3)}_{j} \right\Vert_{L^{2}(\Sigma_{j})} \right)^{2}h^{-40}M^{64-2} \\ \leq CM^{4} \cdot \tilde{C}_{j}\left(\gamma_{j}^{-1} +  e^{\tilde{m}_{j}\gamma_{j}}\left\Vert \partial_{\nu}w^{(1,2,3)}_{j} \right\Vert_{L^{2}(\Sigma_{j})} \right)^{2}h^{-40}M^{62},
\end{split}
\end{equation}
where the factor of $CM^{n+1}$ arises from the measure of
\begin{equation}
\left\{(\tau, \xi_{j}) \in \mathbb{R}^{n + 1}\,|\, \left<\tau, \xi_{j} \right> \leq M \right\} \left(\subset \left\{(\tau, \xi_{j}) \in \mathbb{R}^{n + 1}\,|\, \tau^{2} + |\xi_{j}|^{2} \leq 2M^{2} \right\} \right).
\end{equation}
Since we are in the setting of $\mathbb{R}^{n}\;(n = 2,3)$, we have majorized $M^{n + 1}$ by $M^{4}$. From \eqref{coefficient_smallness_assumption} and Plancherel's Theorem, we have that
\begin{equation} \label{negative_one_norm_step_3}
\vert I_{2} \vert \leq M^{-2}\cdot \left\Vert \chi(\tilde{q}_{j,1} - \tilde{q}_{j,2})\theta^{4}_{h} \right\Vert_{L^{2}(\mathbb{R}^{n + 1})} \leq CM^{-2}
\end{equation}
for a constant $C = C(\Omega, \mathcal{O}^{(j)}, T^{*}, M_{j}) > 0$. Combining \eqref{negative_one_norm_step_2} and \eqref{negative_one_norm_step_3}, we have that
\begin{equation} \label{negative_one_norm_step_1}
\begin{split}
\left \Vert \chi \left(\tilde{q}_{j,1} - \tilde{q}_{j,2} \right) \theta^{4}_{h} \right \Vert_{H^{-1}(\mathbb{R}^{n + 1})}\\ \leq \left(\tilde{C}_{j}h^{-40}M^{66}\left(\gamma_{j}^{-1} +  e^{\tilde{m}_{j}\gamma_{j}}\left\Vert \partial_{\nu}w^{(1,2,3)}_{j} \right\Vert_{L^{2}(\Sigma_{j})} \right)^{2} + CM^{-2}\right)^{1/2}.
\end{split}
\end{equation}
Now, interpolation yields
\begin{equation} \label{interpolation_step_1}
\begin{split}
\left \Vert \left(\tilde{q}_{j,1} - \tilde{q}_{j,2}\right)\theta^{4}_{h} \right \Vert^{2}_{L^{2}(\Omega\setminus\mathcal{O}^{(j)} \times (0, T^{*}))} \leq \left \Vert \chi\left(\tilde{q}_{j,1} - \tilde{q}_{j,2}\right)\theta^{4}_{h} \right \Vert^{2}_{L^{2}(\Omega \times (0, T^{*}))}\\
\leq \left \Vert \chi \left( \tilde{q}_{j,1} - \tilde{q}_{j,2} \right) \theta^{4}_{h} \right \Vert_{H^{-1}(\mathbb{R}^{n+1})}\left \Vert \chi \left( \tilde{q}_{j,1} - \tilde{q}_{j,2} \right) \theta^{4}_{h} \right \Vert_{H^{1}(\mathbb{R}^{n+1})}\\
\leq \tilde{C}_{j}h^{-20}M^{33}\left(\gamma_{j}^{-1} +  e^{\tilde{m}_{j}\gamma_{j}}\left\Vert \partial_{\nu}w^{(1,2,3)}_{j} \right\Vert_{L^{2}(\Sigma_{j})} \right)\cdot h^{-1} + CM^{-1}h^{-1},
\end{split}
\end{equation}
where we used \eqref{coefficient_smallness_assumption} and \eqref{leading_order_amplitude_estimate} in order estimate
\begin{equation}
\left \Vert \chi \left( \tilde{q}_{j,1} - \tilde{q}_{j,2} \right) \theta^{4}_{h} \right \Vert_{H^{1}(\mathbb{R}^{n+1})} \leq Ch^{-1}
\end{equation}
for a constant $C = C(\Omega, \mathcal{O}^{(j)}, T^{*}, M_{j}) > 0$. Thus, we have that
\begin{equation} \label{interpolation_step_2}
\begin{split}
\left \Vert \tilde{q}_{j,1} - \tilde{q}_{j,2}\right \Vert^{2}_{L^{2}(\Omega\setminus\mathcal{O}^{(j)} \times (0, T^{*}))}\\ \leq C\left \Vert \left(\tilde{q}_{j,1} - \tilde{q}_{j,2}\right)\theta^{4}_{h}\right \Vert^{2}_{L^{2}(\Omega\setminus\mathcal{O}^{(j)} \times (0, T^{*}))} + C\left \Vert \left(\tilde{q}_{j,1} - \tilde{q}_{j,2}\right)(1-\theta^{4}_{h})\right \Vert^{2}_{L^{2}(\Omega\setminus\mathcal{O}^{(j)} \times (0, T^{*}))}\\
\leq \tilde{C}_{j}h^{-21}M^{33}\left(\gamma_{j}^{-1} +  e^{\tilde{m}_{j}\gamma_{j}}\left\Vert \partial_{\nu}w^{(1,2,3)}_{j} \right\Vert_{L^{2}(\Sigma_{j})} \right) + CM^{-1}h^{-1} + C_{0}h,
\end{split}
\end{equation}
where we estimated
\begin{equation}
\left \Vert \left(\tilde{q}_{j,1} - \tilde{q}_{j,2}\right)(1-\theta^{4}_{h})\right \Vert^{2}_{L^{2}(\Omega\setminus\mathcal{O}^{(j)} \times (0, T^{*}))} \leq C_{0}h
\end{equation}
for a constant $C_{0} = C_{0}(T^{*}, M_{j}) > 0$. This uses \eqref{coefficient_smallness_assumption} and the fact from \eqref{smooth_cut-off_function_1} that $1-\theta^{4}_{h}$ is supported on $[0, 2h] \cup [T^{*} - 2h, T^{*}]$.\\ \indent
Setting $h \coloneqq M^{-1/2}$ (for $M > 1$ large enough to ensure $h \in (0, \min\{1,T^{*}/4\})$), the estimate in \eqref{interpolation_step_2} becomes
\begin{equation} \label{interpolation_step_3}
\begin{split}
\left \Vert \tilde{q}_{j,1} - \tilde{q}_{j,2}\right \Vert^{2}_{L^{2}(\Omega\setminus\mathcal{O}^{(j)} \times (0, T^{*}))} \\
\leq \tilde{C}_{j}M^{87/2}\left(\gamma_{j}^{-1} +  e^{\tilde{m}_{j}\gamma_{j}}\left\Vert \partial_{\nu}w^{(1,2,3)}_{j} \right\Vert_{L^{2}(\Sigma_{j})} \right) + CM^{-1/2}
\end{split}
\end{equation}
for constants $\tilde{C}_{j} = \tilde{C}_{j}\left(\Omega, \Gamma_{j}, \mathcal{O}^{(j)}, T, T^{*},\Vert p_{j} \Vert_{W^{7,\infty}(Q)} , M_{j} \right) > 0$ and $C = C(\Omega, \mathcal{O}^{(j)}, T^{*}, M_{j}) > 0$. We will also pick $M > 1$ so that it satisfies
\begin{equation} \label{interpolation_step_4}
M^{87/2}\gamma_{j}^{-1} = \gamma_{j}^{-1/2} \left(\iff M \coloneqq \gamma_{j}^{1/87}\right).
\end{equation}
As a result, \eqref{interpolation_step_3} becomes
\begin{equation} \label{interpolation_step_5}
\left \Vert \tilde{q}_{j,1} - \tilde{q}_{j,2} \right \Vert^{2}_{L^{2}(\Omega\setminus\mathcal{O}^{(j)} \times (0, T^{*}))} \leq \tilde{C}_{j}\left(\gamma_{j}^{-1/174} + e^{\tilde{m}_{1,j}\gamma_{j}}\left\Vert \partial_{\nu}w^{(1,2,3)}_{j} \right\Vert_{L^{2}(\Sigma_{j})}\right)
\end{equation}
where $\tilde{m}_{1,j} \geq \tilde{m}_{j}$ is such that $\gamma_{j}^{1/2}e^{\tilde{m}_{j}\gamma_{j}} \leq e^{\tilde{m}_{1,j}\gamma_{j}}$ and $\tilde{C}_{j} = \tilde{C}_{j}\left(\Omega, \Gamma_{j}, \mathcal{O}^{(j)}, T, T^{*},\Vert p_{j} \Vert_{W^{7,\infty}(Q)} , M_{j} \right) > 0$. \\ \indent
We finally consider $\delta_{*,j} \in (0, 1)$ such that when
\begin{equation} \label{linearized_Neumann_data_smallness_0}
\delta_{j,\Gamma_{j}} \coloneqq \left\Vert \partial_{\nu}w^{(1,2,3)}_{j} \right\Vert_{L^{2}(\Sigma_{j})} \in (0, \delta_{*,j}),
\end{equation}
we have that
\begin{equation}
\gamma_{j} \coloneqq \tilde{m}_{1,j}^{-1}\ln \left(\delta_{j,\Gamma_{j}}^{-1/2}\right) \geq \gamma_{*,j} \quad \text{and} \quad h \coloneqq M^{-1/2} = \gamma_{j}^{-1/174} \in (0, \min\{1,T^{*}/4\}).
\end{equation}
Then, \eqref{interpolation_step_5} becomes
\begin{equation} \label{interpolation_step_7}
\left \Vert \tilde{q}_{j,1} - \tilde{q}_{j,2} \right \Vert^{2}_{L^{2}(\Omega\setminus\mathcal{O}^{(j)} \times (0, T^{*}))} \leq \tilde{C}_{j}\left[\left(\ln \left(\delta_{j,\Gamma_{j}}^{-1/2}\right) \right)^{-1/174} + \delta_{j,\Gamma_{j}}^{1/2}\right],
\end{equation}
when \eqref{linearized_Neumann_data_smallness_0} holds. We have $\tilde{C}_{j} = \tilde{C}_{j}\left(\Omega, \Gamma_{j}, \mathcal{O}^{(j)}, T, T^{*},\Vert p_{j} \Vert_{W^{7,\infty}(Q)} , M_{j}\right) > 0$ in \eqref{interpolation_step_7}. \\ \indent
The proof of \eqref{interpolation_step_7} relied on the CGO solutions constructed in \eqref{first_order_linearization_CGO_0}-\eqref{adjoint_CGO_0} for $v^{(l)}_{j}$ and $v_{0,j}$ respectively. In relation to equations \eqref{first_order_linearization_time-dependent} and \eqref{adjoint_problem_CGO_approach}, we shall define 
\begin{equation}
f^{(j)}_{l} \coloneqq v^{(l)}_{j}\big|_{\Sigma},\,f_{0,j} \coloneqq v_{0, j}|_{\Sigma}
\end{equation}
and check that the boundary traces have sufficient regularity to ensure well-posedness by Theorem \ref{linear_well-posedness}. Indeed, we have from Theorem $2.3$ in Chapter $4$ of \cite{LM1972}, equations \eqref{first_order_linearization_CGO_0}-\eqref{adjoint_CGO_0}, the Banach algebra property of $H^{m}(Q)\,(m \geq 3)$, and the estimates in Proposition \ref{LLZ2023_CGO_estimate_proposition} that
\begin{equation} \label{boundary_trace_CGO_estimate_0}
\begin{split}
\left\Vert f^{(j)}_{l} \right\Vert_{H^{9/2,9/2}(\Sigma)} \leq C\left\Vert v^{(l)}_{j} \right\Vert_{H^{5,5}(Q)} \leq C\left\Vert v^{(l)}_{j} \right\Vert_{H^{5}(Q)}\\
\leq C\left(\left\Vert a^{(l)}_{j} \right\Vert_{H^{5}(Q)} + \left\Vert \exp\left\{i\Phi\left(\cdot, \cdot\,; \omega^{(l)}_{j} \right)\right\} \right\Vert_{H^{5}(Q)}\sum_{k' = 1}^{3}\rho^{-k'}\left\Vert a^{(l),j}_{k'} \right\Vert_{H^{5}(Q)} + \left\Vert r^{(l)}_{j} \right\Vert_{H^{5}(Q)}\right)\\
\leq C'\rho^{10}h^{-9} = C'\gamma_{j}^{1749/174}
\end{split}
\end{equation}
where we substituted $\rho = \gamma_{j}, h = \gamma_{j}^{-1/174},$ and $C' = C'\left(\Omega, T, N_{j} \right) > 0$ (recall that $N_{j}$ is defined in \eqref{coefficient_smallness_assumption}). Likewise, we can show that
\begin{equation} \label{boundary_trace_CGO_estimate_1}
\left\Vert f_{0,j} \right\Vert_{H^{9/2,9/2}(\Sigma)} \leq C\rho^{10}\left<\tau, \xi_{j} \right>^{13}h^{-9} = C\left<\tau, \xi_{j} \right>^{13}\gamma_{j}^{1749/174}
\end{equation}
for a constant $C = C(\Omega, T, N_{j}) > 0$. This shows that $\{f^{(j)}_{l}, f_{0,j}\} \subseteq H^{9/2,9/2}(\Sigma)$. The construction of the CGO amplitude and remainder terms in Section \ref{CGO_section} (in particular, recall equations \eqref{smooth_cut-off_function_1_observation} and \eqref{CGO_0_remainder_equation}) ensure that $\{f^{(j)}_{l}, f_{0,j}\} \subseteq \mathcal{H}^{9/2,9/2}(\Sigma)$.\\ \indent
We lastly check that given $\delta_{j,\Gamma_{j}} \in (0, \delta_{*,j})$, we can pick $\epsilon_{1}, \epsilon_{2}, \epsilon_{3} \in \mathbb{R}$ sufficiently small to ensure that
\begin{equation}
\epsilon_{1}f^{(j)}_{1} + \epsilon_{2}f^{(j)}_{2} + \epsilon_{3}f^{(j)}_{3} \in \mathcal{B}_{r}(0),
\end{equation}
where $\mathcal{B}_{r}(0)$ is defined by \eqref{boundary_neighborhood}. By \eqref{boundary_trace_CGO_estimate_0}, it suffices to take
\begin{equation}
\epsilon_{l} \leq \frac{r}{3C'}\gamma_{j}^{-1749/174} = \frac{r}{3C'}\left[\tilde{m}_{1,j}^{-1}\ln \left(\delta^{-1/2}\right)\right]^{-1749/174} \:\,(l = 1,2,3).
\end{equation}
This choice of $\epsilon_{l}$ ensures the well-posedness of \eqref{two_component_GPEqn_small_data} by Theorem \ref{nonlinear_well-posedness}.\\ \indent
The argument in Sections \ref{Fourier_transform_estimate_section}-\ref{completing_the_proof_section} can be repeated in a similar manner for the other nonlinear coefficients in \eqref{two_component_GPEqn_small_data} and this completes the proof of Theorem \ref{main_inverse_stability_theorem}.
\begin{remark} \label{unique_determination_result_0}
By passing derivatives onto test functions in $C^{\infty}(\partial\Omega)$, we can exchange derivatives in the weak sense for the Banach space-valued mapping
\begin{equation*}
t \longmapsto \partial_{\nu}w^{(1,2,3)}_{j}(\cdot, t)\big|_{\partial\Omega} \in H^{3/2}(\partial\Omega), 
\end{equation*}
which justifies \eqref{Neumann_data_measurement_0}. Denoting the solution to \eqref{two_component_GPEqn_small_data} by $u_{j,k}^{(\epsilon_{1}, \epsilon_{2}, \epsilon_{3})}$, this allows one to infer a unique determination result from \eqref{main_stability_estimate_statement} if there exist $\epsilon_{l,0}\;(l = 1,2,3)$ such that
\begin{equation}
\partial_{\nu}u_{j,1}^{(\epsilon_{1}, \epsilon_{2}, \epsilon_{3})}\big|_{\Sigma_{j}} = \partial_{\nu}u_{j,2}^{(\epsilon_{1}, \epsilon_{2}, \epsilon_{3})}\big|_{\Sigma_{j}} 
\end{equation}
for all sufficiently small $\epsilon_{l} \in (-\epsilon_{l,0}, \epsilon_{l,0})\setminus\{0\}\;(l = 1,2,3)$.
\end{remark}
\section{Recovery by Carleman Estimates} \label{Carleman_estimate_recovery_section}
In this section, we consider the setting in which \eqref{two_component_GPEqn_small_data} is equipped with trivial boundary data but nontrivial initial data and outline the main idea for deriving an estimate similar to \eqref{interpolation_step_7} using the method from \cite{AZ2025}. However, the method in \cite{AZ2025} is only able to treat time-independent coefficients and as a result, we assume in what follows that the coefficients $p_{j}, q_{j,k},\tilde{q}_{j,k}$ in \eqref{two_component_GPEqn_small_data} depend only on $x \in \Omega$.
\\ \indent
Denoting $\mathcal{T}_{j} \coloneqq i\partial_{t} + \Delta + p_{j}(x)$, we now consider 
\begin{equation} \label{two_component_GPEqn_small_initial_data}
\begin{cases}
\mathcal{T}_{j} u_{j,k} + q_{j,k}(x)\left|u_{j,k}\right|^{2}u_{j,k} + \tilde{q}_{j,k}(x)|u_{\iota(j),k}|^{2}u_{j,k} = 0 \text{ on } Q,\\
u_{j,k}(x,t) = 0 \text{ on } \Sigma,\\
u_{j}(x,0) = \left(\epsilon_{1}f^{(j)}_{1} + \epsilon_{2}f^{(j)}_{2} + \epsilon_{3}f^{(j)}_{3}\right)(x) \text{ on } \Omega \times \{0\},
\end{cases}
\end{equation}
where $j \in \{1,2\}$, $k \in \{1,2\}$, and $f^{(j)}_{l} \in H^{4}(\Omega) \cap H^{3}_{0}(\Omega)$ for $l \in \{1,2,3\}$. We take $p_{j}, q_{j,k},\tilde{q}_{j,k} \in H^{2}(\Omega)$. We further have the regularity assumptions
\begin{equation} \label{coefficient_smallness_assumption_1}
\begin{split}
\Vert p_{j} \Vert_{H^{2}(\Omega)} \leq M_{0},\\
\Vert q_{j,1} - q_{j,2} \Vert_{H^{4}(\Omega)} \leq M_{1},\\
\Vert \tilde{q}_{j,1} - \tilde{q}_{j,2} \Vert_{H^{4}(\Omega)}\leq M_{2},\\
\left\Vert f^{(j)}_{l}\right\Vert_{H^{4}(\Omega)} \leq \gamma_{+}
\end{split}
\end{equation}
for some constants $M_{0}, M_{1}, M_{2}, \gamma_{+} > 0$. Letting $\mathcal{O}^{(j)} \subset \Omega$ be a relatively open neighborhood of $\partial \Omega$, we also assume that
\begin{equation} \label{coefficient_assumption_0}
\begin{split}
q_{j,1} \equiv q_{j,2} \text{ on } \mathcal{O}^{(j)},\\
\tilde{q}_{j,1} \equiv \tilde{q}_{j,2} \text{ on } \mathcal{O}^{(j)},\\
\left|f^{(j)}_{l}(x)\right| \geq \gamma_{-} \text{ on } \Omega \setminus \mathcal{O}^{(j)}
\end{split}
\end{equation}
for some constant $\gamma_{-} > 0$. Then, we have the following result.
\begin{theorem} \label{main_theorem_other_setting}
Let $\Omega$ be a smooth bounded domain in $\mathbb{R}^{n}\;(n = 2, 3)$. Consider $\Gamma_{j}\;(j = 1,2)$ relatively open subsets of $\partial\Omega$ and let $\mathcal{O}^{(j)} \subset \Omega$ be a relatively open neighborhood of $\partial \Omega$ such that $\partial \mathcal{O}^{(j)} \setminus \partial\Omega$ is $C^{2}$. For \eqref{two_component_GPEqn_small_initial_data}, suppose we have real-valued $f^{(j)}_{l} \in H^{4}(\Omega) \cap H^{3}_{0}(\Omega)$ and the conditions in \eqref{coefficient_smallness_assumption_1}-\eqref{coefficient_assumption_0} satisfied. \\ \indent
Then, for $j \in \{1,2\}$, there exists $\delta_{*,j} > 0$ such that when
\begin{equation}
\delta_{j,\Gamma_{j}} \coloneqq \left\Vert \partial_t\left(\partial^{3}_{\epsilon_{1}\epsilon_{2}\epsilon_{3}}\partial_{\nu}(u_{j,1}-u_{j,2})\Big|_{\epsilon_{1} = \epsilon_{2} = \epsilon_{3} = 0}\right) \right\Vert_{L^{2}(\Gamma_{j} \times (0, T))} \leq \delta_{*,j},
\end{equation}
we have
\begin{equation}
\left \Vert \tilde{q}_{j,1} - \tilde{q}_{j,2} \right \Vert^{2}_{L^{2}(\Omega\setminus\mathcal{O}^{(j)})} \leq C_{j}\left[\left(\ln \left(\delta_{j,\Gamma_{j}}^{-1}\right) \right)^{-1} + \delta_{j,\Gamma_{j}}\right]
\end{equation}
for a constant $C_{j} = C_{j}\left(\Omega, \Gamma_{j}, \mathcal{O}^{(j)}, M_{0}, M_{1}, M_{2}, T, \gamma_{-}, \gamma_{+} \right) > 0$. 
\end{theorem}
\begin{remark} \label{section_5_strategy_technical_remark}
\hfill
\begin{itemize}
\item[(a)] The conditions $p_{j}, q_{j,k}, \tilde{q}_{j,k} \in H^{2}(\Omega)$ and $f^{(j)}_{l} \in H^{4}(\Omega) \cap H^{3}_{0}(\Omega)$ ensure that the solutions to linearized versions of \eqref{two_component_GPEqn_small_initial_data} (namely equations \eqref{section_5_first_order_linearization}-\eqref{partial_t_two_component_GPEqn}) have unique solutions in $C\left([0, T]; H^{2}(\Omega)\cap H^{1}_{0}(\Omega)\right)$. Moreover, the conditions in \eqref{coefficient_smallness_assumption_1} ensure that we can prove estimates similar to those of Lemma $4.2$ in \cite{AZ2025} for the solution to \eqref{partial_t_two_component_GPEqn}.
\item[(b)] The reason we need the conditions in \eqref{coefficient_assumption_0} is related to \eqref{Project_1_approach_0}-\eqref{Project_1_approach_1} in which we need to estimate
\begin{equation} \label{Project_1_approach_0_0}
\left\Vert  \left(\tilde{q}_{j,1} - \tilde{q}_{j,2} \right)\left[f^{(j)}_{3}f^{(\iota(j))}_{1}f^{(\iota(j))}_{2} \right] \right\Vert^{2}_{L^{2}(\Omega_{1})},
\end{equation}
where we recall that $\Omega_{l}\;(l = 1,2)$ was defined in relation to \eqref{boundary_neighborhoods_0}. By the Sobolev Embedding Theorem, $f^{(j)}_{l} \in H^{4}(\Omega) \cap H^{3}_{0}(\Omega)$ becomes arbitrarily small near the boundary. As a result, we need to stipulate that $q_{j,k}$ and $\tilde{q}_{j,k}$ are determined near the boundary. Moreover, converting \eqref{Project_1_approach_0_0} into an estimate of $\left\Vert \tilde{q}_{j,1} - \tilde{q}_{j,2} \right\Vert^{2}_{L^{2}(\Omega_{1})}$ requires a lower bound for the initial data $f^{(j)}_{l}$ away from the boundary.
\end{itemize}
\end{remark}
We will need the following Carleman estimate from \cite{BC2009} (also see Proposition $3.3$ in \cite{AZ2025}):
\begin{proposition} \label{standard_Schrodinger_Carleman_estimate_0}
Let $p(x) \in L^{\infty}(\Omega)$ with $\Vert p\Vert_{L^{\infty}(\Omega)} \leq M$ for some $M > 0$. Picking a point $x_{0} \notin \overline{\Omega}$ and $T^{*} \in (0, T)$, define $\psi(x) = |x - x_{0}|^{2}$ and
\begin{equation}
\phi(x,t) =  \frac{e^{2\lambda\Vert \psi\Vert_{\infty}} - e^{\lambda \psi(x)}}{(T^{*})^{2}-t^{2}}
\end{equation}
on $\Omega \times (-T^{*},T^{*})$, where $\Vert \psi \Vert_{\infty} = \Vert \psi \Vert_{L^{\infty}(\Omega)}$. Now, there exist $\lambda_{*} = \lambda_{*}(\Omega, T^{*}, \Omega_{1}, \Omega_{2}) > 0,\, s_{*} = s_{*}(\Omega, T^{*}, \Omega_{1}, \Omega_{2}) > 0$ such that for all $\lambda \geq \lambda_{*},\, s \geq s_{*},$ and a constant $C = C(M,\Omega, T^{*}, s_{*},\lambda_{*}) > 0$,
\begin{equation} \label{OldCarlemanEstimate_2}
\begin{split}
s^3 \lambda^4\Vert e^{-s\phi}u\Vert_{L^2(\Omega_{2} \times (-T^{*},T^{*}))}^2 + s\lambda\Vert e^{-s\phi}\nabla u\Vert_{L^2(\Omega_{2} \times (-T^{*},T^{*}))}^2
\leq
C\Vert e^{-s\phi} Lu\Vert_{L^2(\Omega \times (-T^{*},T^{*}))}^2 \\ 
+\; C\left(s^3 \lambda^4\Vert e^{-s\phi}u\Vert_{L^2(\Omega_{3} \setminus \Omega_{2} \times (-T^{*},T^{*}))}^2 + s\lambda\Vert e^{-s\phi}\nabla u\Vert_{L^2(\Omega_{3} \setminus \Omega_{2} \times (-T^{*},T^{*}))}^2 \right).
\end{split}
\end{equation}
The above holds for $u \in L^{2}\left(-T^{*},T^{*};\,H^{1}_{0}(\Omega) \right)$ such that $Lu \in L^{2}\left(\Omega \times (-T^{*},T^{*})\right)$.
\end{proposition}
\begin{remark} \label{OldCarlemanEstimate_2_remark}
One can show that \eqref{OldCarlemanEstimate_2} also holds with $\Omega_{3}$ in place of $\Omega_{2}$ on the left-hand side (e.g. see Remark $3.2$ in \cite{AZ2025}). 
\end{remark}
Now, when $\epsilon_{l} \in \mathbb{R}$ is sufficiently small, one can combine the methods from Theorem \ref{nonlinear_well-posedness} and Section $2.2$ in \cite{AZ2025} to show that \eqref{two_component_GPEqn_small_initial_data} has a unique solution in $C\left([0,T]; H^{2}(\Omega) \cap H^{1}_{0}(\Omega)\right)$. Denoting this unique solution to \eqref{two_component_GPEqn_small_initial_data} by $u_{j,k}\left(x,t;\,\epsilon_{1}f^{(j)}_{1} + \epsilon_{2}f^{(j)}_{2} + \epsilon_{3}f^{(j)}_{3}\right)$, one can also show that the map
\begin{equation}
(\epsilon_{1}, \epsilon_{2}, \epsilon_{3}) \in \mathbb{R}^{3} \longmapsto u_{j,k}\left(x,t;\,\epsilon_{1}f^{(j)}_{1} + \epsilon_{2}f^{(j)}_{2} + \epsilon_{3}f^{(j)}_{3}\right)
\end{equation}
is differentiable to third order at $(\epsilon_{1},\epsilon_{2},\epsilon_{3}) = (0, 0, 0)$. So, we can consider
\begin{equation} \label{section_5_first_order_linearization}
v^{(l)}_{j}(x,t) \coloneqq \partial_{\epsilon_{l}}u_{j,k}\big|_{\epsilon_{l} = 0}
\end{equation}
and
\begin{equation} \label{section_5_kth_order_linearization}
w^{(1,2,3)}_{j}(x,t) \coloneqq \partial^{3}_{\epsilon_{1}\epsilon_{2} \epsilon_{3}}\left(u_{j,1} - u_{j,2} \right)\big|_{\epsilon_{1} = \epsilon_{2} = \epsilon_{3} = 0}
\end{equation}
as in Section \ref{rigorous_linearization}. We can show that $r_{j}(x,t) = \partial_{t}w^{(1,2,3)}_{j}$ satisfies
\begin{equation} \label{partial_t_two_component_GPEqn}
\begin{cases}
\mathcal{T}_{j}r_{j} = \partial_{t}\left[\left(\mathcal{S}^{(1,2,3)}_{j,1} - \mathcal{S}^{(1,2,3)}_{j,2} \right) + \left(\tilde{\mathcal{S}}^{(1,2,3)}_{j,1} - \tilde{\mathcal{S}}^{(1,2,3)}_{j,2} \right)\right] \text{ on } Q,\\
r_{j} = 0 \text{ on } \Sigma,\\
r_{j}(x,0) = -i\left[\left(\mathcal{S}^{(1,2,3)}_{j,1} - \mathcal{S}^{(1,2,3)}_{j,2}\right) + \left(\tilde{\mathcal{S}}^{(1,2,3)}_{j,1} - \tilde{\mathcal{S}}^{(1,2,3)}_{j,2} \right)\right](x,0) \text{ on } \Omega \times \{0\},
\end{cases}
\end{equation}
where we recall the notation in \eqref{kth_order_time-dependent_source_terms} for the source terms. 
\begin{remark}
We can define the solution for $t \in (-T, 0)$ by taking $v^{(l)}_{j}(x,t) = \overline{v^{(l)}_{j}(x,-t)}$ (this is why we assumed $f^{(j)}_{l}$ was real-valued in Theorem \ref{main_theorem_other_setting}) and $r_{j}(x,t) = -\overline{r_{j}(x,-t)}$. This allows one to use the Carleman estimate from Proposition \ref{standard_Schrodinger_Carleman_estimate_0} and to bound $\Vert \partial r/\partial \nu\Vert_{L^{2}(\Gamma_{j} \times (-T, T))}$ by a constant times $\Vert \partial r/\partial \nu\Vert_{L^{2}(\Sigma_{j})}$.
\end{remark}
Now, note for example that if we pick initial data where $f^{(j)}_{1} \equiv f^{(j)}_{2} \equiv 0$, we have in \eqref{partial_t_two_component_GPEqn} the initial data
\begin{equation}
r_{j}(x,0) = -2i\left(\tilde{q}_{j,1} - \tilde{q}_{j,2} \right)\left[f^{(j)}_{3}f^{(\iota(j))}_{1}f^{(\iota(j))}_{2} \right]
\end{equation}
and also the fact that the source term can be written in the form $(\tilde{q}_{j,1} - \tilde{q}_{j,2})v_{j}$, where $v_{j}$ depends only on $\left(\tilde{\mathcal{S}}^{(1,2,3)}_{j,1} - \tilde{\mathcal{S}}^{(1,2,3)}_{j,2}\right)$. \\ \indent
We can then proceed as in the calculation after Remark $4.8$ in \cite{AZ2025} to obtain an estimate of $\left \Vert \tilde{q}_{j,1} - \tilde{q}_{j,2} \right \Vert_{L^{2}(\Omega)}$. In the following, we summarize the calculation while omitting the parameter dependence of the constants (further details are provided in \cite{AZ2025}). One can show that $e^{-s\phi}r_{j}(x,t)$ satisfies
\begin{equation}
\begin{split}
\Vert e^{-s\phi(\cdot\,,0)}r_{j}(\cdot\,,0) \Vert^{2}_{L^{2}(\Omega_{2})} \leq 
C\Big(s^{-2}\left\Vert e^{-s\phi(\cdot, 0)} \left(\tilde{q}_{j,1} - \tilde{q}_{j,2} \right) \right\Vert^{2}_{L^{2}(\Omega_{3} \times (-T^{*}, T^{*}))} \left\Vert v_{j} \right\Vert_{L^{\infty}(-T^{*}, T^{*}; H^{2}(\Omega_{3}))}\\ +\; s^{2}\Vert e^{-s\phi}r_{j} \Vert^{2}_{L^{2}(\Omega_{3} \times (-T^{*}, T^{*}))}
+\;  \Vert e^{-s\phi}\nabla r_{j} \Vert^{2}_{L^{2}(\Omega_{3} \times (-T^{*}, T^{*}))}\Big)
\end{split}
\end{equation}
for $s \geq s_{0}$, where $s_{0} > 1$ and $C > 0$ are constants. Using estimates similar to those of Lemma $4.2$ in \cite{AZ2025} to estimate $\left\Vert v_{j} \right\Vert_{L^{\infty}(-T^{*}, T^{*}; H^{2}(\Omega_{2}))}$ and the initial data lower bound in \eqref{coefficient_assumption_0}, we further obtain
\begin{equation} \label{Project_1_approach_0}
\begin{split}
\left\Vert e^{-s\phi(\cdot, 0)} \left(\tilde{q}_{j,1} - \tilde{q}_{j,2} \right) \right\Vert^{2}_{L^{2}(\Omega_{2})} \leq C\Big(s^{-2} \left\Vert e^{-s\phi(\cdot, 0)} \left(\tilde{q}_{j,1} - \tilde{q}_{j,2} \right) \right\Vert^{2}_{L^{2}(\Omega_{3} \times (-T^{*}, T^{*}))}  \\+\; s^{2}\left\Vert e^{-s\phi}r_{j} \right\Vert^{2}_{L^{2}(\Omega_{3} \times (-T^{*}, T^{*}))}
+\;  \Vert e^{-s\phi}\nabla r_{j} \Vert^{2}_{L^{2}(\Omega_{3} \times (-T^{*}, T^{*}))} \Big)
\end{split}
\end{equation}
for $s \geq s_{0}$. Taking $s \geq s_{0}$ large enough to absorb the first term on the right-hand side of \eqref{Project_1_approach_0} into the left-hand side, we now obtain
\begin{equation} \label{Project_1_approach_0_1}
\begin{split}
\left\Vert e^{-s\phi(\cdot, 0)} \left(\tilde{q}_{j,1} - \tilde{q}_{j,2} \right) \right\Vert^{2}_{L^{2}(\Omega_{2})} \leq C\Big(s^{2}\left\Vert e^{-s\phi}r_{j} \right\Vert^{2}_{L^{2}(\Omega_{3} \times (-T^{*}, T^{*}))}
+ \Vert e^{-s\phi}\nabla r_{j} \Vert^{2}_{L^{2}(\Omega_{3} \times (-T^{*}, T^{*}))} \Big).
\end{split}
\end{equation}
Let us note the estimate
\begin{equation} \label{OldCarlemanEstimate2_weight_estimate}
\phi(x,0) \leq \frac{e^{2\lambda \Vert \psi \Vert_{\infty}}}{(T^{*})^{2}}.
\end{equation}
Now, applying \eqref{OldCarlemanEstimate2_weight_estimate}, the Carleman estimate in Proposition \ref{standard_Schrodinger_Carleman_estimate_0} with fixed $\lambda \geq \max\{\lambda_{*}, 1\}$, the unique continuation estimate in Proposition \ref{unique_continuation_result_BC2009}, and the assumption that $\tilde{q}_{j,1} \equiv \tilde{q}_{j,2}$ on $\mathcal{O}^{(j)}$ yield
\begin{equation} \label{Project_1_approach_1}
\begin{split}
\Vert  \tilde{q}_{j,1} - \tilde{q}_{j,2} \Vert^{2}_{L^{2}(\Omega \setminus \mathcal{O}^{(j)})} \leq C\Big(s^{-1} \Vert  \tilde{q}_{j,1} - \tilde{q}_{j,2} \Vert^{2}_{L^{2}(\Omega \setminus \mathcal{O}^{(j)})}\Big) \\
+\;Cs^{2}\left(\gamma_{j}^{-1} + e^{m_{j}\gamma_{j}}\left\Vert \partial r_{j}/\partial \nu\right\Vert_{L^{2}(\Sigma_{j})}^{2}\right)
\end{split}
\end{equation}
for constants $C, \gamma_{*,j}, m_{j} > 0$ and sufficiently large $s \geq s_{0}, \gamma_{j} \geq \gamma_{*,j}$. We can again take $s \geq s_{0}$ large enough to absorb the first term on the right-hand side of \eqref{Project_1_approach_1} into the left-hand side. Now, the required result follows by defining
\begin{equation}
\gamma_{j} = m_{j}^{-1}\ln \left\Vert \partial r_{j}/\partial\nu\right\Vert_{L^{2}(\Sigma_{j})}^{-1}
\end{equation}
for
\begin{equation}
\left\Vert \partial r_{j}/\partial\nu\right\Vert_{L^{2}(\Sigma_{j})} \leq \delta_{*,j},
\end{equation}
which is chosen small enough to ensure that $\gamma_{j} \geq \gamma_{*,j}$. Through considerations similar to those of Remark \ref{unique_determination_result_0}, we can interchange the order of the derivatives in $\partial r_{j}/\partial\nu$ and this completes the proof of Theorem \ref{main_theorem_other_setting}. We recall from Remark \ref{unique_determination_result_0} that the purpose of interchanging the order of the derivatives is to infer a unique determination result. Similar reasoning can be used to obtain estimates for the other nonlinear coefficients in \eqref{two_component_GPEqn_small_initial_data}.

\end{document}